\documentclass[11pt]{article}
\usepackage{graphicx,wrapfig}
\usepackage{flafter}
\usepackage{leftidx}
\usepackage{mathrsfs}
\usepackage{amssymb,amsmath}
\usepackage{bbding}
\usepackage{fancyhdr}
\usepackage{mathrsfs}
\usepackage{color}
\usepackage{stmaryrd}
\usepackage{amsfonts}
\usepackage{latexsym}
\usepackage{psfrag}
\usepackage{graphicx,subfigure}
\usepackage{fancyhdr,graphicx}
\usepackage{multicol}
\usepackage{dsfont}
\usepackage{bbm}
\usepackage{booktabs}
\usepackage[center]{caption2}
\usepackage{cite}
\usepackage{epstopdf}
\usepackage{booktabs}
\usepackage{multirow}
\usepackage{enumerate}
\usepackage{enumitem}

\allowdisplaybreaks
\renewcommand{\leq}{\leqslant}

\date{}
\newcommand{\zd}{\,\mathrm{d}}

\newtheorem{theorem}{Theorem}[section]
\newtheorem{lemma}{Lemma}[section]
\newtheorem{remark}{Remark}
\newtheorem{example}{Example}[section]
\newtheorem{corollary}{Corollary}[section]
\numberwithin{equation}{section}

\newcommand{\abs}[1]{\left|#1\right|}
\newcommand{\absb}[1]{\big|#1\big|}

\newcommand{\bra}[1]{\left(#1\right)}
\newcommand{\brab}[1]{\big(#1\big)}
\newcommand{\braB}[1]{\Big(#1\Big)}

\newcommand{\kbra}[1]{\left[#1\right]}

\newcommand{\mynormb}[1]{\big\|#1\big\|}
\newcommand{\mynormB}[1]{\Big\|#1\Big\|}

\begin{document}
\title{Simple maximum-principle preserving time-stepping methods for time-fractional Allen-Cahn equation}
\author{Bingquan Ji\thanks{Department of Mathematics, Nanjing University of Aeronautics and Astronautics,
211101, P. R. China. Bingquan Ji (jibingquanm@163.com).}
\quad Hong-lin Liao\thanks{Corresponding author. Department of Mathematics,
Nanjing University of Aeronautics and Astronautics,
Nanjing 211106, P. R. China. Hong-lin Liao (liaohl@nuaa.edu.cn, liaohl@csrc.ac.cn)
is supported by a grant 1008-56SYAH18037
from NUAA Scientific Research Starting Fund of Introduced Talent.}
\quad Luming Zhang\thanks{Department of Mathematics, Nanjing University of Aeronautics and Astronautics,
211101, P. R. China. Luming Zhang (zhanglm@nuaa.edu.cn)
is supported by the research grants No. 11571181 from National Natural Science Foundation of China.}}
\date{}
\maketitle
\normalsize

\begin{abstract}
Two fast L1 time-stepping methods, including the backward Euler and stabilized semi-implicit schemes, are suggested for the time-fractional Allen-Cahn equation with Caputo's derivative.
The time mesh is refined near the initial time to resolve the intrinsically initial singularity of solution,
 and unequal time-steps are always incorporated into our approaches so that an adaptive time-stepping strategy
 can be used in long-time simulations.
 It is shown that the proposed schemes using the fast L1 formula preserve the discrete maximum principle.
 Sharp error estimates reflecting the time regularity of solution are established
 by applying the discrete fractional Gr\"{o}nwall inequality and global consistency analysis.
 Numerical experiments are presented to show the effectiveness of our methods and to confirm our analysis.\\


\noindent{\bf AMS subject classiffications.}\;\; 35Q99, 65M06, 65M12, 74A50
\end{abstract}
\section{Introduction}
The phase field models have become popular to describe a host of free-boundary problems in various areas, including material, physical and biology systems \cite{Allen1979A, Cahn1958Free, Clarke1987Origin, Zhao2016A}.
Relevant numerical methods and simulations  are also increasing substantially \cite{Shen2012Numerical,Lee2016Comparison, Kim2017A}. It is well known that the phase field models permit multiple time scales, i.e. an initial dynamics evolves on a fast time scale and later coarsening evolves on a very slow time scale.
It is therefore to consider the adaptive time-stepping strategy \cite{QiaoAn2011,Zhang2013An,Li2017Computationally}, namely, small time steps are utilized
when the energy dissipates rapidly and large time steps are employed otherwise.
These works suggest that
nonuniform time meshes are preferable in the numerical simulations of phase field models.

In comparison with the bright achievement of classical phase field models, in recent years,
there are many researches on building fractional phase field models,
such as time, space and time-space fractional Allen-Cahn equations
{\cite{Inc2018Time,Akagi2016Fractional,Hou2017Numerical, Zheng2017A,Liu2018Time,Zhao2019On}}
to accurately describe anomalous diffusion problems.
Li et al. {\cite{Zheng2017A}} investigated a space-time fractional Allen-Cahn phase-field model
that describes the transport of the fluid mixture of two immiscible fluid phases.
They concluded that the alternative model could provide more accurate description of anomalous diffusion processes
and sharper interfaces than the classical model.
Hou et al. {\cite{Hou2017Numerical}} showed that a fractional in space Allen-Cahn equation could be
viewed an $L^{2}$ gradient flow for the fractional analogue version of Ginzburg-Landau free energy function.
They proved the energy decay property and the maximum principle of continuous problem.
Recently, the authors of \cite{Inc2018Time} considered the symmetry analysis, explicit solution
and convergence analysis of the time-fractional Allen-Cahn and Klein-Gordon equations
with Riemann-Liouville derivative. Zhao et al. \cite{Liu2018Time,Zhao2019On} studied a series of the
time fractional phase field models numerically. The considerable numerical evidences
indicate that the effective free energy of the time fractional phase field models
obeys a similar power law as the integer ones.

The multi-scale nature of time-fractional phase field models prompts us to construct reliable time-stepping methods
on general nonuniform meshes. In this paper, two nonuniform time-stepping schemes
are investigated for the time-fractional Allen-Cahn equation \cite{Zheng2017A,Liu2018Time,Zhao2019On}
\begin{align}
&\partial_{t}^{\alpha}u
=\varepsilon^{2}\Delta{u}-f(u),\quad
\mathbf{x}\in\Omega,\quad{0}<{t}\leq{T},\label{Problem-1}\\
&u(\mathbf{x},0)=u_{0}(\mathbf{x}),\quad \mathbf{x}\in\bar{\Omega},\label{Problem-2}
\end{align}
where $\mathbf{x}=(x,y)^{T}$ and
$\Omega=(a,b)\times(c,d)$ with its closure $\bar{\Omega}$.
The notation $\partial_{t}^{\alpha}:={}_{0}^{C}\!D_{t}^{\alpha}$ in {\eqref{Problem-1}} denotes
the fractional Caputo derivative of order $\alpha$ with respect to $t$,
\begin{align}\label{CaputoDef}
(\partial_{t}^{\alpha}v)(t)
:=(\mathcal{I}_{t}^{1-\alpha}v')(t)=\int_{0}^{t}\omega_{1-\alpha}(t-s)v'(s)\zd{s},\quad 0<\alpha<1,
\end{align}
involving the fractional  Riemann-Liouville integral $\mathcal{I}_{t}^{\mu}$ of order $\mu>0$, that is,
\begin{align}
(\mathcal{I}_{t}^{\mu}v)(t)
:=\int_{0}^{t}\omega_{\mu}(t-s)v(s)\zd{s},\quad\text{where}\quad  \omega_{\mu}(t):=t^{\mu-1}/\Gamma(\mu).
\end{align}
The nonlinear bulk force $f(u)=u^{3}-u$, and the small constant $\varepsilon>{0}$,
called the interaction length, describes
the thickness of the transition boundary between materials. Boundary conditions
are set to be periodic so as not to complicate the analysis with unwanted details.


Very recently, the energy decay laws of time-fractional phase field models,
involving time-fractional Allen-Cahn equation,
time-fractional Cahn-Hilliard equation and time-fractional molecular beam epitaxy models,
are reported in {\cite{Tang2018On}}.
In comparison to the classical physical model, the energy dissipation law of the
time-fractional Allen-Cahn equation {\eqref{Problem-1}} is
\begin{align}\label{FracEnergyDecayLaw}
E(t)\leq{E}(0),
\end{align}
where
\begin{align}\label{FracEnergyExpression}
E(t):=\int_{\Omega}\left[\frac{\varepsilon^{2}}{2}|\nabla{u}|^{2}
+F(u)\right]\zd{\mathbf{x}},\quad
F(u)=\frac{1}{4}(1-u^{2})^{2}.
\end{align}
Also, it possesses a maximum principle, namely,
\begin{align}\label{FracMaximumPrinciple}
|u(\mathbf{x},t)|\leq{1}\;\text{for $t>0$}\quad
\text{if}\quad|u(\mathbf{x},0)|\leq{1}.
\end{align}
To our knowledge, there are few results in the literature on the discrete energy decay law or maximum principle
of numerical approaches for the time-fractional phase field models, especially on nonuniform time meshes.
One of our interests in this paper is to build two nonuniform L1 schemes preserving the maximum principle
of the problem {\eqref{Problem-1}}.

We consider the nonuniform time levels $0=t_{0}<t_{1}<\cdots<t_{k-1}<t_{k}<\cdots<t_{N}=T$ with the time-step sizes $\tau_{k}:=t_{k}-t_{k-1}$ for $1\leq{k}\leq{N}$ and the maximum time-step size $\tau:=\max_{1\leq{k}\leq{N}}\tau_{k}$. Also, let the local time-step ratio $\rho_k:=\tau_k/\tau_{k+1}$ and the maximum step ratio $\rho:=\max_{k\geq 1}\rho_k$. Given a grid function $\{v^{k}\}$, put $\triangledown_{\tau}v^{k}:=v^{k}-v^{k-1}$, $\partial_{\tau}v^{k-\frac12}:=\triangledown_{\tau}v^{k}/\tau_k$ and
$v^{k-\frac{1}{2}}:=(v^{k}+v^{k-1})/2$ for $k\geq{1}$. Always,
let $(\Pi_{1,k}v)(t)$ denote the linear interpolant of a function $v(t)$ at two nodes $t_{k-1}$ and $t_{k}$,
and define a piecewise linear approximation
\begin{align}\label{linear interpolation}
\Pi_{1}v:=\Pi_{1,k}v\quad\text{so that}\quad(\Pi_{1}v)'(t)=\partial_{\tau}v^{k-\frac12}\quad
\text{for $t_{k-1}<{t}\leq t_{k}$ and $k\geq1$}.
\end{align}

\begin{figure}[htb!]
\centering
\includegraphics[width=3.0in]{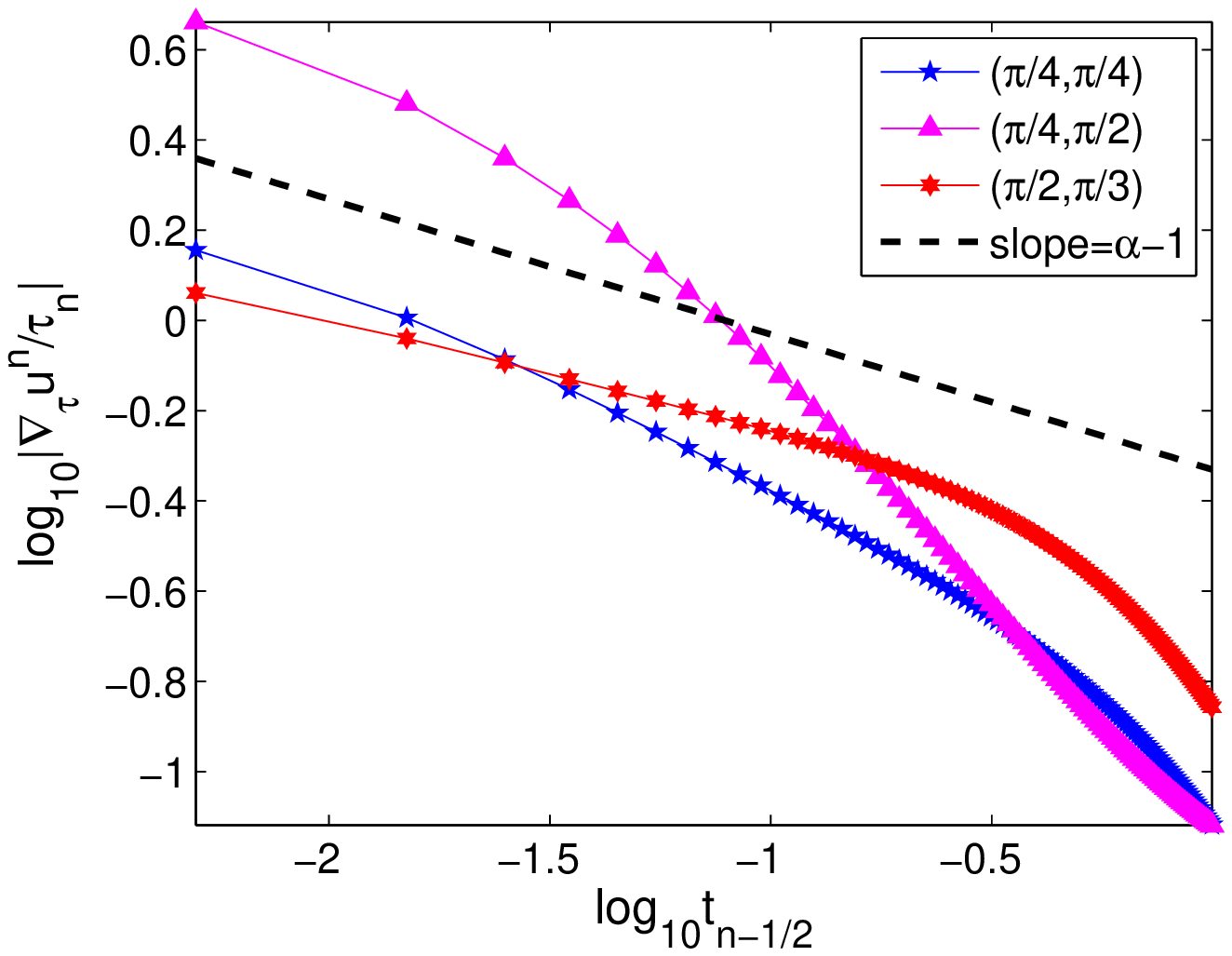}
\includegraphics[width=3.0in]{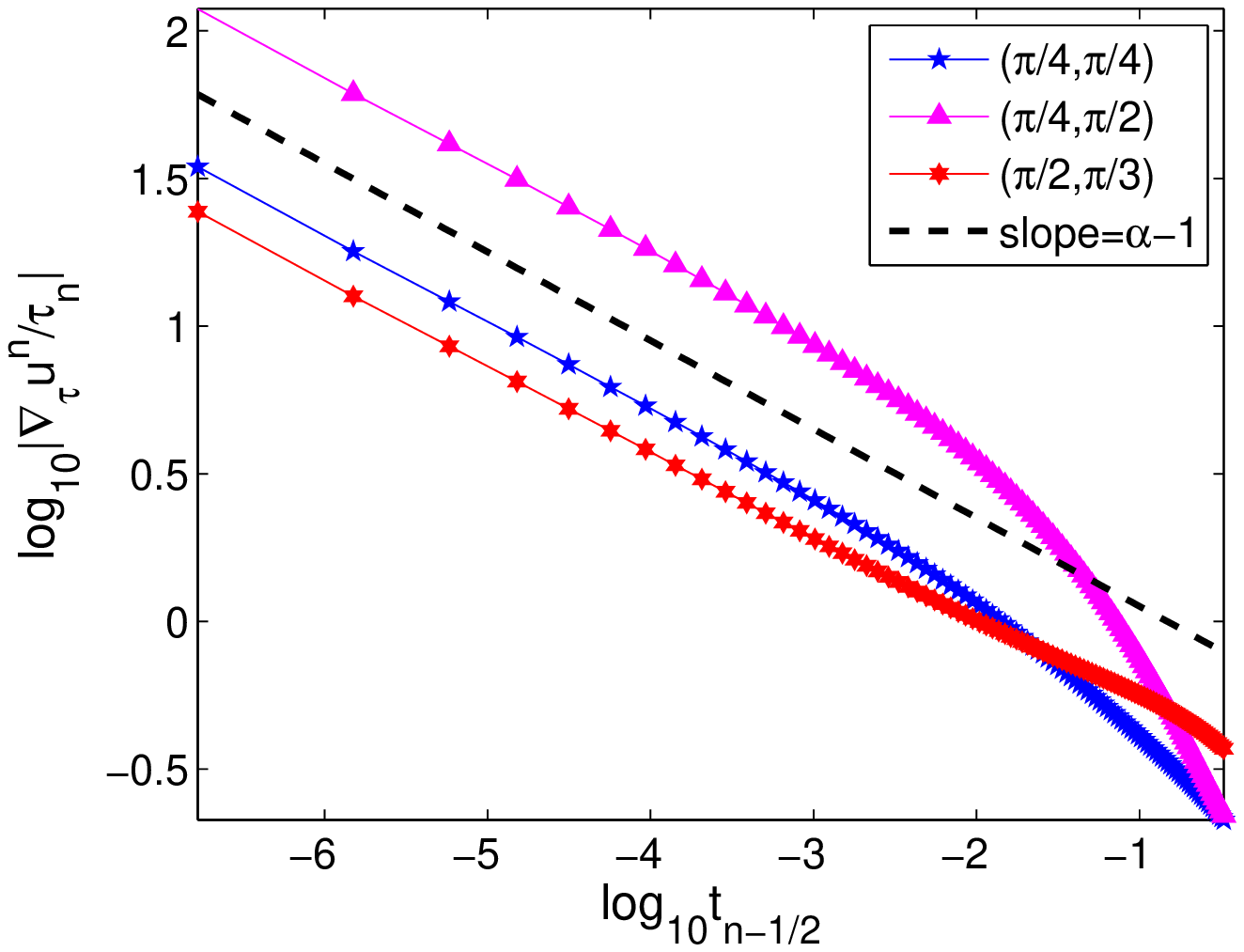}
\caption{The log-log plot of the difference quotient $\partial_{\tau}v^{k-\frac12}$ versus time
for {\eqref{Problem-1}}-{\eqref{Problem-2}} with fractional order $\alpha=0.7$ and $\gamma=1,\,3$ (from left to right), respectively.}
\label{InitialSingularity}
\end{figure}


As an essential mathematical feature of linear and nonlinear subdiffusion problems
including the time-fractional Allen-Cahn problem {\eqref{Problem-1}}-{\eqref{Problem-2}},
the solution always lacks the smoothness near the initial time
although it would be smooth away from  $t=0$, see \cite{Jin2016An,JinLiZhou:2017}.
Actually, assuming the nonlinear function $f$ is Lipschitz continuous
and the initial data $u^0\in H^2(\Omega)\cap H_0^1(\Omega)$,
Jin et al. \cite[Theorem 3.1]{JinLiZhou:2017} proved the subdiffusion problem has an unique solution $u$ for which
$u\in C\bra{[0,T];H^2(\Omega)\cap H_0^1(\Omega)}$, $\partial_t^{\alpha}u\in C\bra{[0,T];L^2(\Omega)}$ and $
\partial_t u \in L^2(\Omega)$  with
$\|\partial_t u(t)\|_{L^2(\Omega)} \le C_ut^{\alpha-1}$ for $0<t\leq T$.
The L1 scheme with a lagging linearized technique for handling the nonlinearity $f(u)$
has been analyzed, and \cite[Theorem 4.5]{JinLiZhou:2017} showed that the discrete solution is $O(\tau^{\alpha})$
convergent in $L^\infty(L^2(\Omega))$. It formally implies that, in any numerical methods for solving
time-fractional diffusion equations, a key consideration is the singularity of the solution near the time $t=0$, see also \cite{Liao2018Sharp,Liao2018Unconditional,Liao2018second}. More directly, we consider the L1 scheme for the time-fractional problem {\eqref{Problem-1}}-{\eqref{Problem-2}}
describing the coalescence of two kissing bubbles, see more details in Example \ref{NumericalApplication}.
Fig. {\ref{InitialSingularity}} plots the discrete time derivative $\partial_{\tau}v^{k-\frac12}$
  near $t=0$ on the graded mesh $t_k=(k/N)^{\gamma}$.
They suggest that
$$\log|u_{t}(\mathbf{x},t)|\approx(\alpha-1)\log{t}+C(\mathbf{x})\quad\text{as $t\rightarrow0$.}$$
It says that the solution possesses weak singularity like $u_{t}=O(t^{\alpha-1})$ near initial time, which can be alleviated by using the graded meshes.
Thus the second interest of this paper is to resolve the essentially weak singularity
in the equation {\eqref{Problem-1}} by refining time mesh near $t=0$. Actually, we will show that
the graded mesh can recover the optimal time accuracy of L1 formula when the solution $u$
does not have the required regularity.

In the next section, we construct the backward Euler and stabilized semi-implicit schemes
by using the nonuniform fast L1 formula $(\partial_{f}^{\alpha}u)^{n}$ described in \eqref{FastApproximation}.
Theorems \ref{DisMaxPrinciple} and \ref{DisMaxPrinciple-stabilized} show that both the backward Euler method \eqref{L1-1}-\eqref{L1-2}
and stabilized semi-implicit method {\eqref{StabilizedL1-1}}-{\eqref{StabilizedL1-2}}
preserve the maximum principle \eqref{FracMaximumPrinciple} in the discrete level
such that they are unconditionally stable in the maximum norm.
By using the recently proposed discrete fractional Gr\"{o}nwall inequality \cite{Liao2018discrete}
and the global consistency analysis \cite{Liao2018Unconditional} of L1 formula,
 we prove that, see Theorems \ref{ConvergenceTheorem} and \ref{ConvergenceTheorem-stabilized}, the fully implicit method
 \eqref{L1-1}-\eqref{L1-2} is convergent
with an optimal order of $O(\tau^{2-\alpha})$ and the stabilized scheme {\eqref{StabilizedL1-1}}-{\eqref{StabilizedL1-2}} is convergent
with an optimal order of $O(\tau)$ in time on the graded meshes with a grading parameter $\gamma\ge1$.
Unfortunately, we are not able to establish any discrete energy dissipation laws on general nonuniform meshes
and leave it as an open problem (see Remark \ref{L1 quadraticForm}).

In summary, the main contributions of this paper are the following:
(i) develop two fast L1 time-stepping methods with unequal time-steps preserving the discrete maximum principle,
(ii) prove the unconditional convergence with the optimal accuracy in time.
Extensive numerical experiments are curried out in section 4 to support our analysis.
Some further remarks conclude the article.

\section{Fast L1 time-stepping methods}
The well-known L1 formula of Caputo derivative \eqref{CaputoDef} is given by
\begin{align}\label{L1-Formula}
(\partial_{\tau}^{\alpha}v)^{n}
:=\int_{t_{0}}^{t_{n}}
\omega_{1-\alpha}(t_{n}-s)(\Pi_{1}v)'(s)\zd{s}
=\sum_{k=1}^{n}a_{n-k}^{(n)}\triangledown_{\tau}v^{k},
\end{align}
where the corresponding discrete convolution kernels $a_{n-k}^{(n)}$ are defined by
\begin{align}\label{L1-Coefficient}
a_{n-k}^{(n)}
:=\frac{1}{\tau_{k}}\int_{t_{k-1}}^{t_{k}}\omega_{1-\alpha}(t_n-s)\zd{s}\quad\text{for $1\leq{k}\leq{n}.$}
\end{align}
Obviously, the discrete convolutional kernels $a_{n-k}^{(n)}$ are positive and decreasing, see also \cite{Liao2018Sharp,Liao2018Unconditional},
\begin{align}\label{L1CoeffEstima}
a_{n-k}^{(n)}>0 \quad\text{and}\quad a_{n-k-1}^{(n)}>a_{n-k}^{(n)}
\quad\text{for $1\leq{k}\leq{n-1}.$}
\end{align}
Note that, this property {\eqref{L1CoeffEstima}} is essential to
the preservation of maximum principle for the proposed L1-type schemes described below.

\subsection{Fast L1 formula}
It is well known that the standard L1 formula \eqref{L1-Formula} is prohibitively expensive for long time simulations. Therefore, to reduce the computational cost and storage requirements incurred by employing the L1 formula directly,
we apply the sum-of-exponentials (SOE) technique to speed up the evaluation of the original problem.
A core result is to approximate the kernel function $t^{-\alpha}$ efficiently on the interval $[\Delta{t},\,T]$,
see \cite[Theorem 2.5]{Jiang2017Fast}.
\begin{lemma}\label{SOE}
For the given $\alpha\in(0,\,1)$, an absolute tolerance error $\epsilon\ll{1}$, a cut-off time $\Delta{t}>0$ and a finial time $T$, there exists a positive integer $N_{q}$, positive quadrature nodes $\theta^{\ell}$ and corresponding positive weights $\varpi^{\ell}\,(1\leq{\ell}\leq{N_{q}})$ such that
\begin{align}
\bigg|
\omega_{1-\alpha}(t)
-\sum_{\ell=1}^{N_{q}}\varpi^{\ell}e^{-\theta^{\ell}t}\bigg|\leq\epsilon,
\quad
\forall\,{t}\in[\Delta{t},\,T].\nonumber
\end{align}
\end{lemma}

To be more precise, the Caputo derivative \eqref{CaputoDef} is split into the sum of a history part (an integral over $[0,\,t_{n-1}]$) and a local part (an integral over $[t_{n-1},\,t_{n}]$) at the time $t_{n}$. Then, the local part will be approximated by linear interpolation directly, the history part can be evaluated via the SOE technique, that is,
\begin{align}\label{FastApprox}
(\partial_{t}^{\alpha}v)(t_{n})
&\approx
\int_{0}^{t_{n-1}}v^{\prime}(s)
\sum_{\ell=1}^{N_{q}}
\varpi^{\ell}e^{-\theta^{\ell}(t_{n}-s)}\zd{s}
+\int_{t_{n-1}}^{t_{n}}\omega_{1-\alpha}(t_{n}-s)
\frac{\triangledown_{\tau}v^{n}}{\tau_{n}}\zd{s}\nonumber\\
&=\sum_{\ell=1}^{N_{q}}\varpi^{\ell}e^{-\theta^{\ell}\tau_{n}}\mathcal{H}^{\ell}(t_{n-1})
+a_{0}^{(n)}\triangledown_{\tau}v^{n},\quad n\geq{1},
\end{align}
where $\mathcal{H}^{\ell}(t_{0}):=0$ and $\mathcal{H}^{\ell}(t_{k}):=\int_{0}^{t_{k}}e^{-\theta^{\ell}(t_{k}-s)}v^{\prime}(s)\zd{s}$.
By utilizing the linear interpolation and a recursive formula, we can approximate $\mathcal{H}^{\ell}(t_{k})$ by
\begin{align}\label{HistoryPart}
\mathcal{H}^{\ell}(t_{k})
&\approx\int_{0}^{t_{k-1}}e^{-\theta^{\ell}(t_{k}-s)}v^{\prime}(s)\zd{s}
+\int_{t_{k-1}}^{t_{k}}e^{-\theta^{\ell}(t_{k}-s)}\frac{\triangledown_{\tau}v^{k}}{\tau_{k}}\zd{s}\nonumber\\
&=e^{-\theta^{\ell}\tau_{k}}\mathcal{H}^{\ell}(t_{k-1})
+b^{(k,l)}\triangledown_{\tau}v^{k},
\end{align}
where the positive coefficients
\[
b^{(k,l)}:=\frac{1}{\tau_{k}}
\int_{t_{k-1}}^{t_{k}}e^{-\theta^{\ell}(t_{k}-s)}\zd{s},\quad k\geq{1}.
\]
Having taken this excursion through {\eqref{FastApprox}}-{\eqref{HistoryPart}}, we arrive at the fast algorithm of L1 formula
\begin{align}\label{FastApproximation}
(\partial_{f}^{\alpha}v)^{n}
=
a_{0}^{(n)}\triangledown_{\tau}v^{n}
+\sum_{\ell=1}^{N_{q}}\varpi^{\ell}e^{-\theta^{\ell}\tau_{n}}\mathcal{H}^{\ell}(t_{n-1}),
\quad n\geq{1},
\end{align}
in which $\mathcal{H}^{\ell}(t_{k})$ is computed by using the recursive relationship
\begin{align}\label{HistoryRecursive}
\mathcal{H}^{\ell}(t_{k})
=e^{-\theta^{\ell}\tau_{k}}\mathcal{H}^{\ell}(t_{k-1})
+b^{(k,l)}\triangledown_{\tau}v^{k},\quad k\geq{1},\quad 1\leq\ell\leq{N_{q}}.
\end{align}

For the convenience of numerical analysis, we now eliminate the historic term $\mathcal{H}^{\ell}(t_{k})$
from the fast L1
formula \eqref{FastApproximation}. From the recursive equation \eqref{HistoryRecursive},
a direct calculation yields
\begin{align}\label{HistoryIdentity}
\mathcal{H}^{\ell}(t_{k})
=\sum_{j=1}^{k}e^{-\theta^{\ell}(t_{k}-t_{j})}b^{(j,l)}\triangledown_{\tau}v^{j},\quad k\geq{1},\quad 1\leq{l}\leq{N_{q}}.
\end{align}
By substituting \eqref{HistoryIdentity} into {\eqref{FastApproximation}}, we get an alternative definition
\begin{align}\label{FastL1-Formula}
(\partial_{f}^{\alpha}v)^{n}
=a_{0}^{(n)}\triangledown_{\tau}v^{n}
+\sum_{k=1}^{n-1}\frac{\triangledown_{\tau}v^{k}}{\tau_{k}}
\int_{t_{k-1}}^{t_{k}}\sum_{\ell=1}^{N_{q}}\varpi^{\ell}e^{-\theta^{\ell}(t_{n}-s)}\zd{s}
=\sum_{k=1}^{n}A_{n-k}^{(n)}\triangledown_{\tau}v^{k},\quad n\geq{1},
\end{align}
where the corresponding discrete convolution coefficient $A_{n-k}^{(n)}$ is defined by
\begin{align}\label{FastL1-Coefficient}
A_{0}^{(n)}:=a_{0}^{(n)},\quad
A_{n-k}^{(n)}:=\frac{1}{\tau_{k}}
\int_{t_{k-1}}^{t_{k}}\sum_{\ell=1}^{N_{q}}\varpi^{\ell}e^{-\theta^{\ell}(t_{n}-s)}\zd{s},\quad
{1}\leq{k}\leq{n-1},\,{n}\geq{1}.
\end{align}
For the discrete kernels $A_{j}^{(n)}$,
we have the following result \cite[Lemma 2.5]{Liao2018Unconditional}.
\begin{lemma}\label{lemma:FastL1CoeffEstimate}
If the tolerance error $\epsilon$ of SOE satisfies $\epsilon\leq\min{\{\frac{1}{3}\omega_{1-\alpha}(T),\,\alpha\omega_{2-\alpha}(1)\}}$, then the discrete convolutional kernel $A_{n-k}^{(n)}$ of {\eqref{FastL1-Coefficient}} satisfies
\begin{itemize}
  \item [(i)] $A_{k-1}^{(n)}>A_{k}^{(n)}>0$ for $1\leq{k}\leq{n-1}$;
 \item [(ii)] $A_{0}^{(n)}=a_{0}^{(n)}$ and \,$A_{n-k}^{(n)}\geq\frac{2}{3}a_{n-k}^{(n)}$ for $1\leq{k}\leq{n-1}$.
\end{itemize}
\end{lemma}

\subsection{Backward Euler scheme}

We recall briefly the  difference approximation in space.
For two positive integers $M_{1},\,M_{2}$,
let the spatial lengths $h_{1}:=(b-a)/M_{1},\,h_{2}:=(d-c)/M_{2}$ and
 $x_{i}=a+ih_{1}$, $y_{j}=c+jh_{2}$ for $0\leq i\leq M_1$, $0\leq j\leq M_2$.
Also, denote $\bar{\Omega}_{h}:=\big\{\mathbf{x}_{h}=(x_{i},y_{j})\,|\,0\leq i\leq M_1, 0\leq j\leq M_2\}$
and put $\Omega_{h}:=\bar{\Omega}_{h}\cap\Omega$.
For any grid function $\{v_h\,|\,\mathbf{x}_{h}\in\bar{\Omega}_{h}\}$, denote a grid function space
\[
\mathbb{V}_{h}:=\big\{v\,|\,v=(v_{j})^{T}\;\;\text{for}\;\;0\leq{j}\leq{M_{2}-1},
\;\;\text{with}\;\;
v_{j}=(v_{i,j})^{T}\;\;\text{for}\;\;0\leq{i}\leq{M_{1}-1}\big\},
\]
where $v^{T}$ is the transpose of the vector $v$.
The maximum norm $\|v\|_{\infty}:=\max_{\mathbf{x}_{h}\in\Omega_{h}}|v_{h}|$.

Let $D_{h}$ be the discrete matrix of Laplace operator $\Delta$ subject to periodic boundary conditions.
With the Kronecker tensor product $\otimes$, the matrix $D_{h}=I_{1}\otimes{D_{1}}+D_{2}\otimes{I_{2}}$,
in which $I_{1}$ and $I_{2}$ are the identity matrices of order $M_{2}\times{M}_{2}$ and $M_{1}\times{M}_{1}$, respectively,
 and the matrices $D_{1}$ and $D_{2}$ are of forms
\[
D_{1}
=\frac{1}{h_{1}^{2}}
\left(
\begin{array}{ccccc}
-2 &  1 & 0 & \cdots & 1 \\
1  & -2 & 1 & \cdots & 0 \\
\vdots & \ddots & \ddots & \ddots & \vdots\\
0  & \cdots & 1 & -2 & 1 \\
1  & \cdots & 0 & 1  &-2 \\
\end{array}
\right)_{M_{1}\times{M}_{1},}\,
D_{2}
=\frac{1}{h_{2}^{2}}
\left(
\begin{array}{ccccc}
-2 &  1 & 0 & \cdots & 1 \\
1  & -2 & 1 & \cdots & 0 \\
\vdots & \ddots & \ddots & \ddots & \vdots\\
0  & \cdots & 1 & -2 & 1 \\
1  & \cdots & 0 & 1  &-2 \\
\end{array}
\right)_{M_{2}\times{M}_{2}.}
\]
Then we have some primary properties of the discrete matrix $D_{h}$ in the next lemma,
which is straightforward to check and we thus omit the proof here.
\begin{lemma}\label{NegativeCondition}
Under the periodic boundary condition, the discrete matrix $D_{h}$ of the Laplace operator possesses the following properties
\begin{itemize}
  \item [(a)] The discrete matrix $D_{h}$ is symmetric.
  \item [(b)] For any nonzero $v\in{\mathbb{V}_{h}}$, $v^{T}D_{h}v\leq{0}$, i.e., the matrix $D_{h}$ is negative semi-definite.
  \item [(c)] The elements of $D_{h}=(d_{ij})$ fulfill $d_{ii}=-\max_{i}\sum_{j\neq{i}}|d_{ij}|$ for each $i$.
\end{itemize}
\end{lemma}

Now we have the backward Euler-type scheme on irregular meshes
for {\eqref{Problem-1}}-{\eqref{Problem-2}},
\begin{align}
\brab{\partial_{f}^{\alpha}u}^{n}
&=\varepsilon^{2}D_{h}u^{n}
-f(u^{n}),\quad
{n}\geq{1},\label{L1-1}\\
u_{h}^{0}&=u_{0}(\mathbf{x}_{h}),\quad
\mathbf{x}_{h}\in\bar{\Omega}_{h},\label{L1-2}
\end{align}
where $f(u^{n}):=(u^{n})^{3}-u^{n}$ with the vector $(u^{n})^{3}=\bra{(u_{1}^{n})^{3},(u_{2}^{n})^{3},\cdots,(u_{M_{2}-1}^{n})^{3}}^{T}$ and
\[
(u_{j}^{n})^{3}
=\bra{(u_{1,j}^{n})^{3},(u_{2,j}^{n})^{3},\cdots,(u_{M_{1}-1,j}^{n})^{3}}^{T}\quad
\text{for $j=0,1,\cdots,M_2-1$.}
\]

Now we prove that the fully discrete scheme {\eqref{L1-1}}-{\eqref{L1-2}} preserves the maximum principle numerically.
Always, we need the following result \cite[Lemma 3.2]{Hou2017Numerical}.
\begin{lemma}\label{MatrixInfNorm}
Let $B$ be a real $M\times{M}$ matrix and $A=aI-B$ with $a>0$.
If the elements of $B=(b_{ij})$ fulfill $b_{ii}=-\max_{i}\sum_{j\neq{i}}|b_{ij}|$,
then for any $c>0$ and $V\in{\mathbb{R}^{M}}$ we have
\begin{align}
\|AV\|_{\infty}\geq{a}\|V\|_{\infty},\quad
\|AV+c(V)^{3}\|_{\infty}\geq{a}\|V\|_{\infty}+c\|V\|_{\infty}^{3}.\nonumber
\end{align}
\end{lemma}

\begin{theorem}\label{DisMaxPrinciple}
If $\|u^{0}\|_{\infty}\leq{1}$ and the maximum time-step size $\tau\leq1/\sqrt[\alpha]{\Gamma(2-\alpha)},$
then the solution of backward Euler scheme {\eqref{L1-1}}-{\eqref{L1-2}} satisfies $\|u^{k}\|_{\infty}\leq{1}$
for $0\leq{k}\leq{N}.$
So it preserves the maximum principle \eqref{FracMaximumPrinciple} numerically
and is unconditionally stable.
\end{theorem}
\begin{proof}
We use the complete mathematical induction.
Obviously, the claimed inequality holds for $k=0$.
For $1\leq n\le N$, assume that
\begin{align}\label{inductionAssume}
\|u^{k}\|_{\infty}\leq{1}\quad\text{for $0\leq{k}\leq{n-1}.$}
\end{align}
It remains to verify that $\|u^{n}\|_{\infty}\leq{1}$. From the definition \eqref{FastL1-Formula}, one has
\begin{align}\label{L1-Formula-alternative}
(\partial_{f}^{\alpha}u)^{n}
=A_{0}^{(n)}u^{n}-L^{n-1}\quad\text{where}\quad L^{n-1}:=\sum_{k=1}^{n-1}\big(A_{n-k-1}^{(n)}-A_{n-k}^{(n)}\big)u^{k}+A_{n-1}^{(n)}u^{0}.
\end{align}
Thanks to the decreasing property in Lemma \ref{lemma:FastL1CoeffEstimate} (i), the induction hypothesis \eqref{inductionAssume}
and the triangle inequality yield
\begin{align*}
\big\|L^{n-1}\big\|_{\infty}\leq \sum_{k=1}^{n-1}\big(A_{n-k-1}^{(n)}-A_{n-k}^{(n)}\big)\|u^{k}\|_{\infty}+A_{n-1}^{(n)}\|u^{0}\|_{\infty}
\leq A_{0}^{(n)} .
\end{align*}
Then, from the numerical scheme {\eqref{L1-1}}, it is easy to obtain
\begin{align}\label{mid Estimate}
\big\|(A_{0}^{(n)}-1)u^{n}+(u^{n})^{3}
-\varepsilon^{2}D_{h}u^{n}\big\|_{\infty}=\big\|L^{n-1}\big\|_{\infty}\leq A_{0}^{(n)}.
\end{align}
For the left hand side of \eqref{mid Estimate}, we apply Lemma \ref{NegativeCondition} (c) and Lemma {\ref{MatrixInfNorm}}
to find that
\begin{align*}
\big\|(A_{0}^{(n)}-1)u^{n}+(u^{n})^{3}
-\varepsilon^{2}D_{h}u^{n}\big\|_{\infty}
\geq
(A_{0}^{(n)}-1)\big\|u^{n}\big\|_{\infty}
+\big\|u^{n}\big\|_{\infty}^{3}.
\end{align*}
Then it follows from \eqref{mid Estimate} that $(A_{0}^{(n)}-1)\big\|u^{n}\big\|_{\infty}
+\big\|u^{n}\big\|_{\infty}^{3}\leq A_{0}^{(n)}.$
If $A_{0}^{(n)}\ge1$ or the maximum step size $\tau\leq1/\sqrt[\alpha]{\Gamma(2-\alpha)},$
the above inequality implies $\|u^{n}\|_{\infty}\leq1$ immediately.
Otherwise, we have $(A_{0}^{(n)}-1)\|u^{n}\|_{\infty}
+\|u^{n}\|_{\infty}^{3}>A_{0}^{(n)},$
because the function
$$g(z):=(A_{0}^{(n)}-1)z
+z^{3}-A_{0}^{(n)}\quad\text{for $z>0$}$$
is monotonically increasing for any $z>0$. This leads to a contradiction and then
the claimed result holds for $k=n$. The principle of induction completes the proof.
\end{proof}

\subsection{Stabilized semi-implicit scheme}

The backward Euler scheme {\eqref{L1-1}}-{\eqref{L1-2}} is a fully nonlinear implicit scheme and
some inner iteration will be needed. To accelerate the time-stepping process, we build a linearized scheme here by using
the well-known stabilized technique
via a stabilized term $S(u^{n}-u^{n-1})$ for a properly large scalar parameter $S>0$
, see also the recent work \cite{Tang2018On}.
The resulting stabilized semi-implicit scheme for the problem {\eqref{Problem-1}}-{\eqref{Problem-2}} reads
\begin{align}
\brab{\partial_{f}^{\alpha}u}^{n}
&=\varepsilon^{2}D_hu^{n}
-f(u^{n-1})-S(u^{n}-u^{n-1}),\quad
{n}\geq{1},\label{StabilizedL1-1}\\
u_h^0&=u_{0}(\mathbf{x}_h),\quad
\mathbf{x}_h\in\bar{\Omega}_h.\label{StabilizedL1-2}
\end{align}

We have the following result on discrete maximum principle and stability.
\begin{theorem}\label{DisMaxPrinciple-stabilized}
If $\|u^{0}\|_{\infty}\leq{1}$ and the scalar stabilized parameter $S\ge2$,
then the solution of stabilized semi-implicit scheme {\eqref{StabilizedL1-1}}-{\eqref{StabilizedL1-2}} satisfies
\begin{align*}
\|u^{k}\|_{\infty}\leq{1}\quad\text{for $0\leq{k}\leq{N}.$}
\end{align*}
So it preserves the maximum principle \eqref{FracMaximumPrinciple} numerically
and is unconditionally stable.
\end{theorem}
\begin{proof}
It only needs to verify that $\|u^{n}\|_{\infty}\leq{1}$ under the induction hypothesis
\begin{align*}
\|u^{k}\|_{\infty}\leq{1}\quad\text{for $0\leq{k}\leq{n-1}.$}
\end{align*}
From the linearized scheme {\eqref{StabilizedL1-1}}, one has
\begin{align}\label{Tang-mid Estimate}
\big\|(A_{0}^{(n)}+S)u^{n}-\varepsilon^{2}D_{h}u^{n}\big\|_{\infty}=\big\|(1+S)u^{n-1}-(u^{n-1})^{3}+L^{n-1}\big\|_{\infty},
\end{align}
where $L^{n-1}$ is defined in \eqref{L1-Formula-alternative}.
Thanks to the decreasing property in Lemma \ref{lemma:FastL1CoeffEstimate} (i),
the induction hypothesis
and the triangle inequality yield $\big\|L^{n-1}\big\|_{\infty}\leq A_{0}^{(n)}.$
Furthermore, it is easy to check that 
\begin{align*}
\abs{(1+S)z-z^{3}}\leq S\quad\text{if $\abs{z}\leq1$ and $S\ge2$,}
\end{align*}
thus the right hand side of \eqref{Tang-mid Estimate} can be bounded by
\begin{align*}
\big\|(1+S)u^{n-1}-(u^{n-1})^{3}+L^{n-1}\big\|_{\infty}\leq A_{0}^{(n)}+
\big\|(1+S)u^{n-1}-(u^{n-1})^{3}\big\|_{\infty}\le A_{0}^{(n)}+S.
\end{align*}
For the left hand side of \eqref{Tang-mid Estimate}, we apply Lemma \ref{NegativeCondition} (c) and Lemma {\ref{MatrixInfNorm}}
to find that
\begin{align*}
\big\|(A_{0}^{(n)}+S)u^{n}
-\varepsilon^{2}D_{h}u^{n}\big\|_{\infty}
\geq
(A_{0}^{(n)}+S)\big\|u^{n}\big\|_{\infty}.
\end{align*}
Then the desired estimate $\|u^{n}\|_{\infty}\leq{1}$ follows from \eqref{Tang-mid Estimate} directly.
\end{proof}

Due to the presence of the stabilized term $S(u^{n}-u^{n-1})$, the numerical solution generated by
the semi-implicit scheme
{\eqref{StabilizedL1-1}}-{\eqref{StabilizedL1-2}} will be
limited to first-order accurate in time even if the solution is sufficiently smooth.
We address the error analysis in the next section.

\section{Global consistency analysis and convergence}
To facilitate the error analysis of difference approximations in space,  we assume that the continuous solution $u$
is sufficiently smooth in space and satisfies
\begin{align}\label{eq: Regularity sigma}
\|u(t)\|_{W^{4,\infty}(\Omega)}\le C_u,\;\; \|u^{(\ell)}(t)\|_{W^{0,\infty}(\Omega)}\le C_u\brab{1+t^{\sigma-\ell}}\quad \text{for $0<t\leq T$ and $\ell=1,2$,}
\end{align}
where a regularity parameter $\sigma\in(0,1)$ is introduced to make our analysis extendable.

In \cite{Liao2018Unconditional}, the local consistency error
$\Upsilon^{j}:=(\partial_{t}^{\alpha}u)(t_{j})-(\partial_{f}^{\alpha}u)^{j}$
of fast L1 formula \eqref{FastL1-Formula} was bounded by a discrete convolution structure,
which is valid for any time meshes. It provides us an opportunity to give the global error
via the global consistency error $\sum_{j=1}^np_{n-j}^{(n)}\absb{\Upsilon^j}$,
where $p_{n-j}^{(n)}$ are the discrete complementary convolution kernels
defined via \eqref{eq: discreteConvolutionKernel-RL}.
Note that, the definition \eqref{FastL1-Coefficient} and Lemma \ref{lemma:FastL1CoeffEstimate} (i)
show that the discrete convolutional kernels $A_{n-k}^{(n)}$ fulfill two assumptions
\textbf{Ass1}-\textbf{Ass2} in Appendix \ref{sec:gronwall} with $\pi_a=\frac{3}{2}$. In this section,
we will use the results
of Lemma \ref{FractGronwall} without further declarations.

\begin{lemma}\label{lemma:L1formulaNonuniform-consistence}
Under the condition of Lemma \ref{lemma:FastL1CoeffEstimate},
the global consistency error is bounded by
\begin{align*}
\sum_{j=1}^np_{n-j}^{(n)}\absb{\Upsilon^j}
\leq\sum_{k=1}^{n}p_{n-k}^{(n)}A_{0}^{(k)}G^k
+\sum_{k=1}^{n-1}p_{n-k}^{(n)}A_{0}^{(k)}G^k
+\frac{C_{u}}{\sigma}t_{n}^{\alpha}\hat{t}_{n-1}^{2}\epsilon\quad \text{for $n\geq1$,}
\end{align*}
where the local quantities $G^k :=2\int_{t_{k-1}}^{t_k}\bra{t-t_{k-1}}\abs{u_{tt}}\zd t$ for $1\leq k\leq n$ and $\hat{t}_{n}:=\max\{1,t_{n}\}$.
\end{lemma}
\begin{proof} On the basis of the upper bound of $(\partial_{t}^{\alpha})u(t_{n})-(\partial_{\tau}^{\alpha}u)^{n}$ given in \cite[Lemma 3.1]{Liao2018Unconditional},
  the estimate (3.5) in the proof of \cite[Lemma 3.3]{Liao2018Unconditional} gives the desired result.
\end{proof}

To resolve such a solution $u$ efficiently, it is appropriate to choose the time
mesh such that the following
condition~\cite{Liao2018Sharp,Liao2018Unconditional,Liao2018second,William2007A} holds.
\begin{enumerate}[itemindent=1em]
\item[\textbf{AssG}.] Let $\gamma\geq{1}$ be a user-chosen parameter. There is a mesh-independent constant $C_{\gamma}>0$ such that
$\tau_k\le C_{\gamma}\tau\min\{1,t_k^{1-1/\gamma}\}$
for~$1\le k\le N$~and
$t_{k}\leq C_{\gamma}t_{k-1}$ for $2\leq{k}\leq{N}$.
\end{enumerate}
Here, the parameter $\gamma\ge1$ controls the extent to which the time levels
are concentrated near $t=0$. If the mesh is quasi-uniform, then
\textbf{AssG} holds with~$\gamma=1$.  As~$\gamma$ increases, the initial step sizes
become smaller compared to the later ones.
A simple example of a family of meshes satisfying \textbf{AssG} is the graded
mesh~$t_k=T(k/N)^{\gamma}$ with the maximum step ratio $\rho=1$.

It is to note that, the global consistency error in
Lemma \ref{lemma:L1formulaNonuniform-consistence}
gives a super\/convergence estimate of nonuniform L1 formula.
Consider the first time level $n=1$, the regularity setting
\eqref{eq: Regularity sigma} gives
$\abs{\Upsilon^{1}}\leq C_u A_{0}^{(1)}\int_{0}^{t_1}t^{\sigma-1}\zd t\leq
C_u\tau_1^{\sigma-\alpha}/\sigma,$
implying that the L1 formula is always inconsistent if $0<\sigma\le\alpha$,
also see Table \ref{L1-InitialSingularity-1} in Section 4.
However, we have the global consistency error of order $O(\tau_1^{\sigma})$, because
$p_{0}^{(1)}\abs{\Upsilon^1}\le G^1\le C_u\tau_1^{\sigma}/\sigma$.
In general, we have the following result from \cite[Lemma 3.3]{Liao2018Unconditional}.
\begin{corollary}\label{GlobalConsisError}
Under the regularity \eqref{eq: Regularity sigma}, the global consistency error can be bounded by
\begin{align*}
\sum_{j=1}^{n}p_{n-j}^{(n)}\absb{\Upsilon^{j}}
\leq C_u\braB{\,\frac{\tau_1^{\sigma}}{\sigma}
+\frac{1}{1-\alpha}\max_{2\leq{k}\leq{n}}t_{k}^{\alpha}t_{k-1}^{\sigma-2}\tau_{k}^{2-\alpha}
+\frac{\epsilon}{\sigma}t_{n}^{\alpha}\hat{t}_{n-1}^{2}}
\quad\text{for $1\leq{n}\leq{N}.$}
\end{align*}
Specifically, if the mesh satisfies \emph{\textbf{AssG}}, then
\begin{align*}
\sum_{j=1}^{n}p_{n-j}^{(n)}\absb{\Upsilon^{j}}
\leq \frac{C_u}{\sigma(1-\alpha)}\tau^{\min\{2-\alpha,\gamma\sigma\}}
+C_u\frac{\epsilon}{\sigma}t_{n}^{\alpha}\hat{t}_{n-1}^{2}\quad\text{for $1\leq{n}\leq{N}.$}
\end{align*}
\end{corollary}

\begin{theorem}\label{ConvergenceTheorem}
Assume that $\|u_{0}\|_{L^{\infty}(\Omega)}\leq{1}$ and
the solution of \eqref{Problem-1}-\eqref{Problem-2} satisfies the regular assumption \eqref{eq: Regularity sigma}.
If the maximum step size $\tau\leq{1}/\sqrt[\alpha]{6\Gamma(2-\alpha)}$, then the numerical solution $u_{h}^{n}$
of the backward Euler scheme \eqref{L1-1}-\eqref{L1-2} is convergent in the maximum norm, that is,
\begin{align*}
\|u(\mathbf{x}_h,t_n)-u_{h}^{n}\|_{\infty}\leq C_u\braB{\frac{\tau_1^{\sigma}}{\sigma}
+\frac{1}{1-\alpha}\max_{2\leq{k}\leq{n}}t_{k}^{\alpha}t_{k-1}^{\sigma-2}\tau_{k}^{2-\alpha}
+\frac{\epsilon}{\sigma}t_{n}^{\alpha}\hat{t}_{n-1}^{2}
+h_1^{2}+h_2^{2}}
\end{align*}
for $1\leq{n}\leq{N}.$ Moreover, when the time mesh satisfies \emph{\textbf{AssG}}, it holds that
\begin{align*}
\|u(\mathbf{x}_h,t_n)-u_{h}^{n}\|_{\infty}
\leq \frac{C_u}{\sigma(1-\alpha)}
\bra{\tau^{\min\{2-\alpha,\gamma\sigma\}}+\epsilon}
+C_u\bra{h_1^{2}+h_2^{2}}\quad\text{for $1\leq{n}\leq{N},$}
\end{align*}
which achieves the optimal accuracy $O(\tau^{2-\alpha})$ if the graded parameter $\gamma\geq\max{\{1,\,(2-\alpha)/\sigma\}}$.
\end{theorem}
\begin{proof} Let $U_{h}^{n}:=u(\mathbf{x}_h,t_n)$ and the error function $e_{h}^{n}:=U_{h}^{n}-u_{h}^{n}\in{\mathbb{V}_{h}}$
for $\mathbf{x}_{h}\in\bar{\Omega}_{h}$ and $0\leq{n}\leq{N}$.
It is easy to find that the exact solution $U_{h}^{n}$ satisfies the governing equations
\begin{align*}
\bra{\partial_{f}^{\alpha}U}^{n}-\varepsilon^{2}D_{h}U^{n}&=-f(U^{n})+\bra{R_{t}}^{n}+\bra{R_{s}}^{n},\quad {1}\leq{n}\leq{N},\\
U_{h}^{0}&=u_{0}(\mathbf{x}_{h}),\quad \mathbf{x}_{h}\in\Omega_{h},
\end{align*}
where $\bra{R_{t}}^{n}$ and $\bra{R_{s}}^{n}$ denote the truncation errors in time and space, respectively.
Subtracting \eqref{L1-1}-\eqref{L1-2} from the above two equations, respectively,  one gets
\begin{align}
\bra{\partial_{f}^{\alpha}e}^{n}-\varepsilon^{2}D_{h}e^{n}&=
-f(U^{n})+f(u^{n})+\bra{R_{t}}^{n}+\bra{R_{s}}^{n},\quad{1}\leq{n}\leq{N},\label{ErrorEquation-1}\\
e_{h}^{0}&=0,\quad\mathbf{x}_{h}\in\Omega_{h}.\label{ErrorEquation-2}
\end{align}
Recalling the elementary inequality $|(a^{3}-a)-(b^{3}-b)|\leq{2}|a-b|$ for $\forall\, a,b\in[-1,1]$,
we apply Theorem \ref{DisMaxPrinciple} (discrete maximum principle) to get
\begin{align*}
\mynormb{f(U^{n})-f(u^{n})}_{\infty}\leq 2\mynormb{e^{n}}_{\infty}.
\end{align*}
Thus the triangle inequality with the error equation \eqref{ErrorEquation-1} gives
\begin{align}\label{FirstSecondEstimate}
\mynormb{\bra{\partial_{f}^{\alpha}e}^{n}-\varepsilon^{2}D_{h}e^{n}}_{\infty}
\leq 2\mynormb{e^{n}}_{\infty}+\mynormb{\bra{R_{t}}^{n}}_{\infty}+\mynormb{\bra{R_{s}}^{n}}_{\infty}.
\end{align}
Applying the decreasing property $(i)$ of the kernels $A_{n-k}^{(n)}$ and the triangle inequality,
we can bound the left hand side of \eqref{FirstSecondEstimate} by
\begin{align*}
\mynormb{\bra{\partial_{f}^{\alpha}e}^{n}-\varepsilon^{2}D_{h}e^{n}}_{\infty}
&=\mynormB{(A_{0}^{(n)}-\varepsilon^{2}D_{h})e^{n}
-\sum_{k=1}^{n-1}\brab{A_{n-k-1}^{(n)}-A_{n-k}^{(n)}}e^{k}-A_{0}^{(n)}e^{0}}_{\infty}\nonumber\\
&\geq\mynormb{(A_{0}^{(n)}-\varepsilon^{2}D_{h})e^{n}}_{\infty}
-\sum_{k=1}^{n-1}\brab{A_{n-k-1}^{(n)}-A_{n-k}^{(n)}}\mynormb{e^{k}}_{\infty}
-A_{n-1}^{(n)}\mynormb{e^{0}}_{\infty}\nonumber\\
&\geq\sum_{k=1}^{n}A_{n-k}^{(n)}\triangledown_{\tau}\mynormb{e^{k}}_{\infty},
\end{align*}
where Lemma \ref{NegativeCondition} (c) and Lemma {\ref{MatrixInfNorm}} have been used.
Then it follows from \eqref{FirstSecondEstimate} that
\begin{align*}
\sum_{k=1}^{n}A_{n-k}^{(n)}\triangledown_{\tau}\mynormb{e^{k}}_{\infty}
\leq 2\mynormb{e^{n}}_{\infty}+\mynormb{\bra{R_{t}}^{n}}_{\infty}+\mynormb{\bra{R_{s}}^{n}}_{\infty},
\end{align*}
which takes the form of \eqref{eq: first Gronwall} with the substitutions $\lambda:=2$, $v^k:=\mynormb{e^{k}}_{\infty}$,
$\xi^{n}:=\mynormb{\bra{R_{t}}^{n}}_{\infty}$ and $\eta^n:=\mynormb{\bra{R_{s}}^{n}}_{\infty}.$
Lemma \ref{FractGronwall} (the discrete fractional Gr\"{o}nwall inequality) says that, if the maximum step size $\tau\leq{1}/\sqrt[\alpha]{6\Gamma(2-\alpha)}$, then it holds that
\begin{align*}
\mynormb{e^{n}}_{\infty}\leq
2E_{\alpha}\brab{6\max(1,\rho){t}_{n}^{\alpha}}\Big(\max_{1\leq{k}\leq{n}}\sum_{j=1}^{k}p_{k-j}^{(k)}\mynormb{\bra{R_{t}}^{j}}_{\infty}
+\omega_{1+\alpha}(t_{n})\max_{1\leq{k}\leq{n}}\mynormb{\bra{R_{s}}^{k}}_{\infty}\Big).
\end{align*}
Then Corollary \ref{GlobalConsisError} yields the claimed estimate and completes the proof.
\end{proof}

For the semi-implicit scheme \eqref{StabilizedL1-1}-\eqref{StabilizedL1-2}, the global error is dominated
by the stabilized term $S(u^k-u^{k-1})$. Under the regular assumption \eqref{eq: Regularity sigma},
the local consistency error is about $\int_{t_{k-1}}^{t_k}\abs{u_t}\zd{t}$. One can follow the proof of
\cite[Lemma 3.3]{Liao2018Unconditional} to bound
the corresponding global error as follows (also see the case of $m=0$ in the estimate \eqref{eq: P bound})
$$\sum_{j=1}^{n}p_{n-j}^{(n)}\int_{t_{j-1}}^{t_j}\abs{u_t}\zd{t}\leq C_u\braB{\frac{\tau_1^{\sigma}}{\sigma}
+\frac{1}{1-\alpha}\max_{2\leq{k}\leq{n}}t_{k}^{\alpha}t_{k-1}^{\sigma-1}\tau_{k}}.$$
Then, a similar proof of Theorem \ref{ConvergenceTheorem} leads to the following result.
\begin{theorem}\label{ConvergenceTheorem-stabilized}
Assume that $\|u_{0}\|_{L^{\infty}(\Omega)}\leq{1}$ and
the exact solution of \eqref{Problem-1}-\eqref{Problem-2} satisfies the regular assumption \eqref{eq: Regularity sigma}.
If the stabilized parameter $S\ge2$ and the maximum time-step size $\tau\leq{1}/\sqrt[\alpha]{6\Gamma(2-\alpha)}$, then
the numerical solution $u_{h}^{n}$
of the semi-implicit scheme \eqref{StabilizedL1-1}-\eqref{StabilizedL1-2} is convergent in the maximum norm, that is,
\begin{align*}
\|u(\mathbf{x}_h,t_n)-u_{h}^{n}\|_{\infty}\leq C_u\braB{\frac{\tau_1^{\sigma}}{\sigma}
+\frac{1}{1-\alpha}\max_{2\leq{k}\leq{n}}t_{k}^{\alpha}t_{k-1}^{\sigma-1}\tau_{k}
+\frac{\epsilon}{\sigma}t_{n}^{\alpha}\hat{t}_{n-1}^{2}
+h_1^{2}+h_2^{2}}
\end{align*}
for $1\leq{n}\leq{N}$. Moreover, when the time mesh satisfies \emph{\textbf{AssG}}, it holds that
\begin{align*}
\|u(\mathbf{x}_h,t_n)-u_{h}^{n}\|_{\infty}\leq \frac{C_u}{\sigma(1-\alpha)}
\braB{\tau^{\min\{1,\gamma\sigma\}}+\epsilon}
+C_u\bra{h_1^{2}+h_2^{2}}\quad\text{for $1\leq{n}\leq{N},$}
\end{align*}
which achieves the optimal accuracy $O(\tau)$ if the graded parameter $\gamma\geq\max{\{1,\,1/\sigma\}}$.
\end{theorem}

\begin{remark}(An open problem)\label{L1 quadraticForm}
It is interesting to mention that, on the uniform mesh, the discrete L1 kernels \eqref{L1-Coefficient} reads
\begin{align*}
a_{n-k}^{(n)}=a_{n-k}=\frac{1}{\tau^{\alpha}}\kbra{\omega_{2-\alpha}(n-k+1)-\omega_{2-\alpha}(n-k)}\quad\text{for $1\leq{k}\leq{n},$}
\end{align*}
the semi-implicit stabilized scheme \eqref{StabilizedL1-1}-\eqref{StabilizedL1-2} using the L1 formula inherits a discrete energy dissipation law,
see  \cite[Theorem 3.1]{Tang2018On} for details.
As seen, the proof of discrete energy dissipation law relies on
the property of a quadratic form $\sum_{k=1}^nw_k\sum_{j=1}^ka_{k-j}w_j\ge0$.
However, it seems rather difficult to extend the positive semi-definite property to a general class of nonuniform meshes.
More precisely, we are not able to verify the positive semi-definite property of the following quadratic form
(by taking $w_k=\triangledown_{\tau}v^k$)
\begin{align}\label{positive quadratic form}
\sum_{k=1}^n\triangledown_{\tau}v^{k}(\partial_{\tau}^{\alpha}v)^{k}
=\sum_{k=1}^nw_k\sum_{j=1}^ka_{k-j}^{(k)}w_j\ge0.
\end{align}
More generally, it has yet to be determined what restrictions must be imposed on the discrete convolution coefficients $\{A_{n-k}^{(n)}\,|1\leq k\leq n\}$ so that the quadratic form $\sum_{k=1}^nw_k\sum_{j=1}^kA_{k-j}^{(k)}w_j$ is positive semi-definite.
This problem could be challenging and remains open to us.
\end{remark}


\section{Numerical examples}

The nonuniform fast L1 time-stepping methods \eqref{L1-1}-\eqref{L1-2} and
\eqref{StabilizedL1-1}-\eqref{StabilizedL1-2} are examined
for solving the Allen-Cahn problem \eqref{Problem-1}-\eqref{Problem-2}.
Always, we set the absolute tolerance error $\epsilon=10^{-12}$ for the SOE approximation.
The second-order centered difference scheme is used to approximate
the Laplace operator with the same length $h=1/M$ in each spatial direction. For the nonlinear scheme
\eqref{L1-1}-\eqref{L1-2}, a simple iteration
is employed to solve the nonlinear algebra equations at each time level with the termination error $10^{-12}$.
The maximum norm error $e(M,N):=\max_{1\leq{n}\leq{N}}\|U^{n}-u^{n}\|_{\infty}$ is recorded in each run,
and the experimental convergence order in time is computed by
$$\text{Order}:=\frac{\log\bra{e(M,N)/e(M,2N)}}{\log\bra{\tau(N)/\tau(2N)}}$$
where $\tau(N)$ denotes the maximum time-step size for total $N$ subintervals.

\begin{example}\label{ConvergenceTest}
To examine the temporal accuracy of our time-stepping schemes, consider the time-fractional Allen-Cahn equation
$\partial_{t}^{\alpha}u=\frac1{8\pi^{2}}\Delta{u}-f(u)+g(\mathbf{x},t)$ for $\mathbf{x}\in(0,1)^{2}$ and $0<t<1$
such that it has an exact solution $u=\omega_{1+\sigma}(t)\sin(2\pi{x})\sin(2\pi{y})$.
\end{example}

The time interval $[0,T]$ is always divided into two parts $[0, T_0]$ and $[T_0, T]$ with total $N$ subintervals.
We will take $T_0=\min\{1/\gamma,T\}$, and apply the graded grid $t_{k}=T_{0}(k/N_0)^{\gamma}$
in $[0,T_{0}]$ to resolve the initial singularity. In the remainder interval $[T_{0},T]$,
we put $N_1:=N-N_0$ cells with random time-steps
$$\tau_{N_{0}+k}=\frac{(T-T_{0})\epsilon_{k}}{\sum_{k=1}^{N_1}\epsilon_{k}}\quad\text{for $1\leq k\leq N_1$}$$
where $\epsilon_{k}\in(0,1)$ are the random numbers.

\begin{table}[htb!]
\begin{center}
\caption{Temporal error of \eqref{L1-1}-\eqref{L1-2} for $\alpha=0.8,\,\sigma=0.8$ with $\gamma_{\mathrm{opt}}=1.5$}\label{L1-InitialSingularity-1} \vspace*{0.3pt}
\def\temptablewidth{1.0\textwidth}
{\rule{\temptablewidth}{0.5pt}}
\begin{tabular*}{\temptablewidth}{@{\extracolsep{\fill}}cccccccccc}
\multirow{2}{*}{$N$} &\multirow{2}{*}{$\tau$} &\multicolumn{2}{c}{$\gamma=1.25$} &\multirow{2}{*}{$\tau$} &\multicolumn{2}{c}{$\gamma=1.5$} &\multirow{2}{*}{$\tau$}&\multicolumn{2}{c}{$\gamma=2$} \\
             \cline{3-4}          \cline{6-7}         \cline{9-10}
         &          &$e(N)$   &Order &         &$e(N)$   &Order &         &$e(N)$    &Order\\
\midrule
  64     &2.60e-02  &3.57e-03 &$-$   &2.54e-02 &2.65e-03 &$-$   &2.98e-02 &2.33e-03  &$-$\\
  128    &1.25e-02  &1.83e-03 &0.91  &1.32e-02 &1.24e-03 &1.15  &1.42e-02 &9.79e-04	 &1.07\\
  256    &6.44e-03  &9.18e-04 &1.04  &6.76e-03 &5.68e-04 &1.17  &7.10e-03 &4.32e-04	 &1.18\\
  512    &3.15e-03  &4.59e-04 &0.97  &3.46e-03 &2.59e-04 &1.17  &3.61e-03 &1.94e-04	 &1.19\\
\midrule
\multicolumn{3}{l}{$\min\{\gamma\sigma,2-\alpha\}$}   &1.00 & & &1.20 & & &1.20\\
\end{tabular*}
{\rule{\temptablewidth}{0.5pt}}
\end{center}
\end{table}	
\begin{table}[htb!]
\begin{center}
\caption{Temporal error of \eqref{L1-1}-\eqref{L1-2} for $\alpha=0.8,\,\sigma=0.4$ with $\gamma_{\mathrm{opt}}=3$}\label{L1-InitialSingularity-2} \vspace*{0.3pt}
\def\temptablewidth{1.0\textwidth}
{\rule{\temptablewidth}{0.5pt}}
\begin{tabular*}{\temptablewidth}{@{\extracolsep{\fill}}cccccccccc}
\multirow{2}{*}{$N$} &\multirow{2}{*}{$\tau$} &\multicolumn{2}{c}{$\gamma=2$} &\multirow{2}{*}{$\tau$} &\multicolumn{2}{c}{$\gamma=3$} &\multirow{2}{*}{$\tau$}&\multicolumn{2}{c}{$\gamma=4$} \\
             \cline{3-4}          \cline{6-7}         \cline{9-10}
         &          &$e(N)$   &Order &         &$e(N)$   &Order &         &$e(N)$    &Order\\
\midrule
  64     &2.85e-02  &2.67e-02 &$-$   &2.81e-02 &1.75e-02 &$-$   &2.65e-02 &2.13e-02  &$-$\\
  128    &1.45e-02  &1.55e-02 &0.82  &1.36e-02 &8.38e-03 &1.02  &1.40e-02 &1.01e-02	 &1.17\\
  256    &7.22e-03  &8.96e-03 &0.79  &7.23e-03 &3.86e-03 &1.22  &6.83e-03 &4.63e-03	 &1.09\\
  512    &3.68e-03  &5.17e-03 &0.82  &3.66e-03 &1.73e-03 &1.18  &3.51e-03 &2.01e-03	 &1.25\\
\midrule
\multicolumn{3}{l}{$\min\{\gamma\sigma,2-\alpha\}$}   &0.80 & & &1.20 & & &1.20\\
\end{tabular*}
{\rule{\temptablewidth}{0.5pt}}
\end{center}
\end{table}	
\begin{table}[htb!]
\begin{center}
\caption{Temporal error of \eqref{StabilizedL1-1}-\eqref{StabilizedL1-2} for $\alpha=0.8,\,\sigma=0.8$ with $\gamma_{\mathrm{opt}}=1.25$}\label{StabilizedL1-InitialSingularity-1} \vspace*{0.3pt}
\def\temptablewidth{1.0\textwidth}
{\rule{\temptablewidth}{0.5pt}}
\begin{tabular*}{\temptablewidth}{@{\extracolsep{\fill}}cccccccccc}
\multirow{2}{*}{$N$} &\multirow{2}{*}{$\tau$} &\multicolumn{2}{c}{$\gamma=1$} &\multirow{2}{*}{$\tau$} &\multicolumn{2}{c}{$\gamma=1.25$} &\multirow{2}{*}{$\tau$}&\multicolumn{2}{c}{$\gamma=2$} \\
             \cline{3-4}          \cline{6-7}         \cline{9-10}
         &          &$e(N)$    &Order &         &$e(N)$   &Order &         &$e(N)$    &Order\\
\midrule
  64     &1.56e-02	 &1.26e-02 &$-$   &2.87e-02	&9.16e-03 &$-$   &3.70e-02 &7.90e-03  &$-$\\
  128    &7.81e-03	 &6.49e-03 &0.95  &1.47e-02	&4.59e-03 &1.03  &1.84e-02 &3.84e-03  &1.03\\
  256    &3.91e-03	 &3.33e-03 &0.96  &7.69e-03	&2.26e-03 &1.09	 &8.97e-03 &1.88e-03  &0.99\\
  512    &1.95e-03	 &1.70e-03 &0.97  &3.55e-03	&1.11e-03 &0.92  &4.33e-03 &9.19e-04  &0.98\\
\midrule
\multicolumn{3}{l}{$\min\{\gamma\sigma,1\}$}   &0.80 & & &1.00 & & &1.00\\
\end{tabular*}
{\rule{\temptablewidth}{0.5pt}}
\end{center}
\end{table}	
\begin{table}[htb!]
\begin{center}
\caption{Temporal error of \eqref{StabilizedL1-1}-\eqref{StabilizedL1-2} for $\alpha=0.8,\,\sigma=0.4$ with $\gamma_{\mathrm{opt}}=2.5$}\label{StabilizedL1-InitialSingularity-2} \vspace*{0.3pt}
\def\temptablewidth{1.0\textwidth}
{\rule{\temptablewidth}{0.5pt}}
\begin{tabular*}{\temptablewidth}{@{\extracolsep{\fill}}cccccccccc}
\multirow{2}{*}{$N$} &\multirow{2}{*}{$\tau$} &\multicolumn{2}{c}{$\gamma=2$} &\multirow{2}{*}{$\tau$} &\multicolumn{2}{c}{$\gamma=2.5$} &\multirow{2}{*}{$\tau$}&\multicolumn{2}{c}{$\gamma=3$} \\
             \cline{3-4}          \cline{6-7}         \cline{9-10}
         &          &$e(N)$   &Order &         &$e(N)$   &Order &         &$e(N)$    &Order\\
\midrule
  64     &3.74e-02	&2.42e-02 &$-$   &3.55e-02 &1.69e-02 &$-$   &4.00e-02 &1.45e-02  &$-$\\
  128    &1.76e-02  &1.37e-02 &0.75  &1.77e-02 &8.04e-03 &1.06  &1.86e-02 &6.77e-03	 &1.00\\
  256    &8.50e-03  &7.90e-03 &0.76  &9.20e-03 &3.88e-03 &1.12  &9.57e-03 &3.09e-03	 &1.18\\
  512    &4.50e-03  &4.53e-03 &0.87  &4.61e-03 &1.94e-03 &1.01  &4.85e-03 &1.40e-03	 &1.16\\
\midrule
\multicolumn{3}{l}{$\min\{\gamma\sigma,1\}$}   &0.80 & & &1.00 & & &1.00\\
\end{tabular*}
{\rule{\temptablewidth}{0.5pt}}
\end{center}
\end{table}	

We take the spatial grid points $M=1024$ in each direction such that the temporal error dominates the spatial error in each run.
Numerical results of the backward Euler scheme \eqref{L1-1}-\eqref{L1-2}
for two different cases $\sigma=\alpha$ and $\sigma<\alpha$ are listed in Tables \ref{L1-InitialSingularity-1}-\ref{L1-InitialSingularity-2}, respectively. They suggest the time accuracy is of order
  $O(\tau^{\min\{\gamma\sigma,2-\alpha\}})$ and confirm Theorem \ref{ConvergenceTheorem} experimentally.
We also run the stabilized semi-implicit scheme \eqref{StabilizedL1-1}-\eqref{StabilizedL1-2} by setting a variety of regularity parameters.
Tables \ref{StabilizedL1-InitialSingularity-1}-\ref{StabilizedL1-InitialSingularity-2} report
the numerical results in the case $\sigma=\alpha$ and a worse case of $\sigma<\alpha$.
It seen that it is accurate of order $O(\tau^{\min\{\gamma\sigma,1\}})$ on the graded meshes,
confirming Theorem \ref{ConvergenceTheorem-stabilized} experimentally.

\begin{example}[Coalescence of two kissing bubbles]\label{NumericalApplication}
Consider the time-fractional Allen-Cahn problem \eqref{Problem-1}-\eqref{Problem-2}
describing the coalescence of two kissing bubbles inside the spatial domain $\Omega=(-\pi,\,\pi)^{2}$,
by taking $\varepsilon=0.1$ and the initial data
\begin{align*}
u_{0}(\mathbf{x})
=\begin{cases}
0.5, &(x+1)^{2}+y^{2}\leq{1}\; \mathrm{or}\;(x-1)^{2}+y^{2}\leq{1},\\
-0.5, &\mathrm{otherwise}.
\end{cases}
\end{align*}
\end{example}


This example is used to examine the physical effect of the fractional order $\alpha$ in the original problem
and the physical property of our suggested methods.
Theorems \ref{DisMaxPrinciple} and \ref{DisMaxPrinciple-stabilized} suggest that variable time-steps are always allowed in our time-stepping approaches. As a matter of fact, the temporal evolution of phase models involve multiple time scales which initial data evolves on a fast time scale at the early stage of dynamics and then the coarsening evolves rather slowly until it reaches a steady state. Hence, to capture the fast dynamics and reduce the cost of computation, we adapt the variant adaptive time-stepping strategy \cite{Li2017Computationally}
\begin{align*}
\tau_{k}=\min\left\{\max\Big\{\tau_{\min},\,
\frac{tol}{1+\beta\|u^{k}-u^{k-1}\|_{\infty}}\Big\},\tau_{\max}\right\}\quad \text{for $k\geq{1}$},
\end{align*}
where the constant 1 is set to avoid the possible singularity as the model reaches the steady state.
The parameters $tol$ and $\beta$ are used to adjust the level of adaptively and would be chosen in experience.
A small $tol$ or a large $\beta$ will generate time steps close to $\tau_{\min}$, which a large $tol$ or a small $\beta$ will give time steps close to $\tau_{\max}$.
The problem is simulated to the final time $T=100$ by taking $M=128,T_{0}=0.1$, $t_{k}=T_{0}(k/N_{0})^{\gamma}$ with the graded parameter $\gamma=3$ in the initial interval $[0,T_{0}]$ and adopting adaptive time steps in the remainder interval.

\begin{figure}[htb!]
\centering
\includegraphics[width=2.0in]{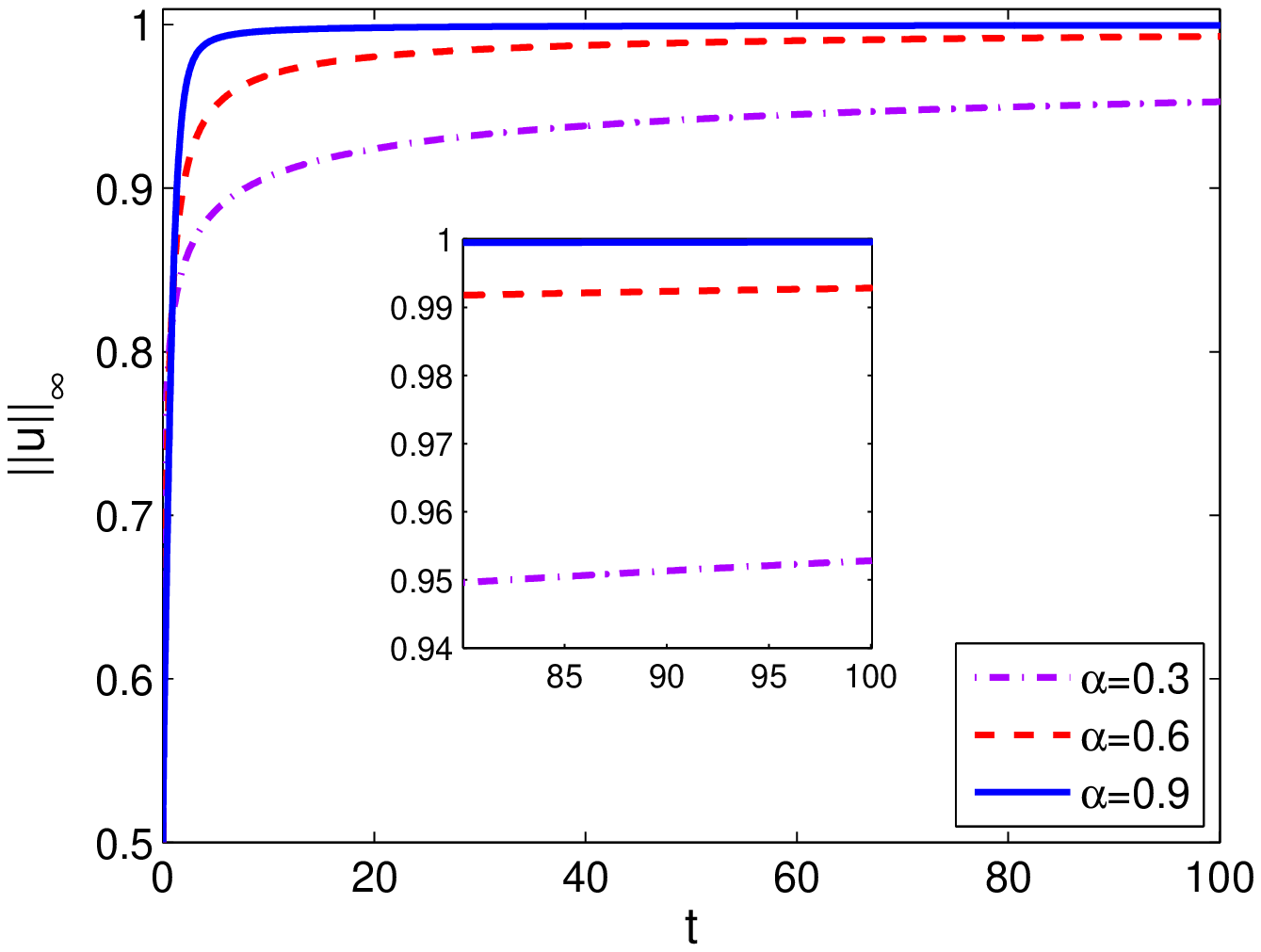}
\includegraphics[width=2.0in]{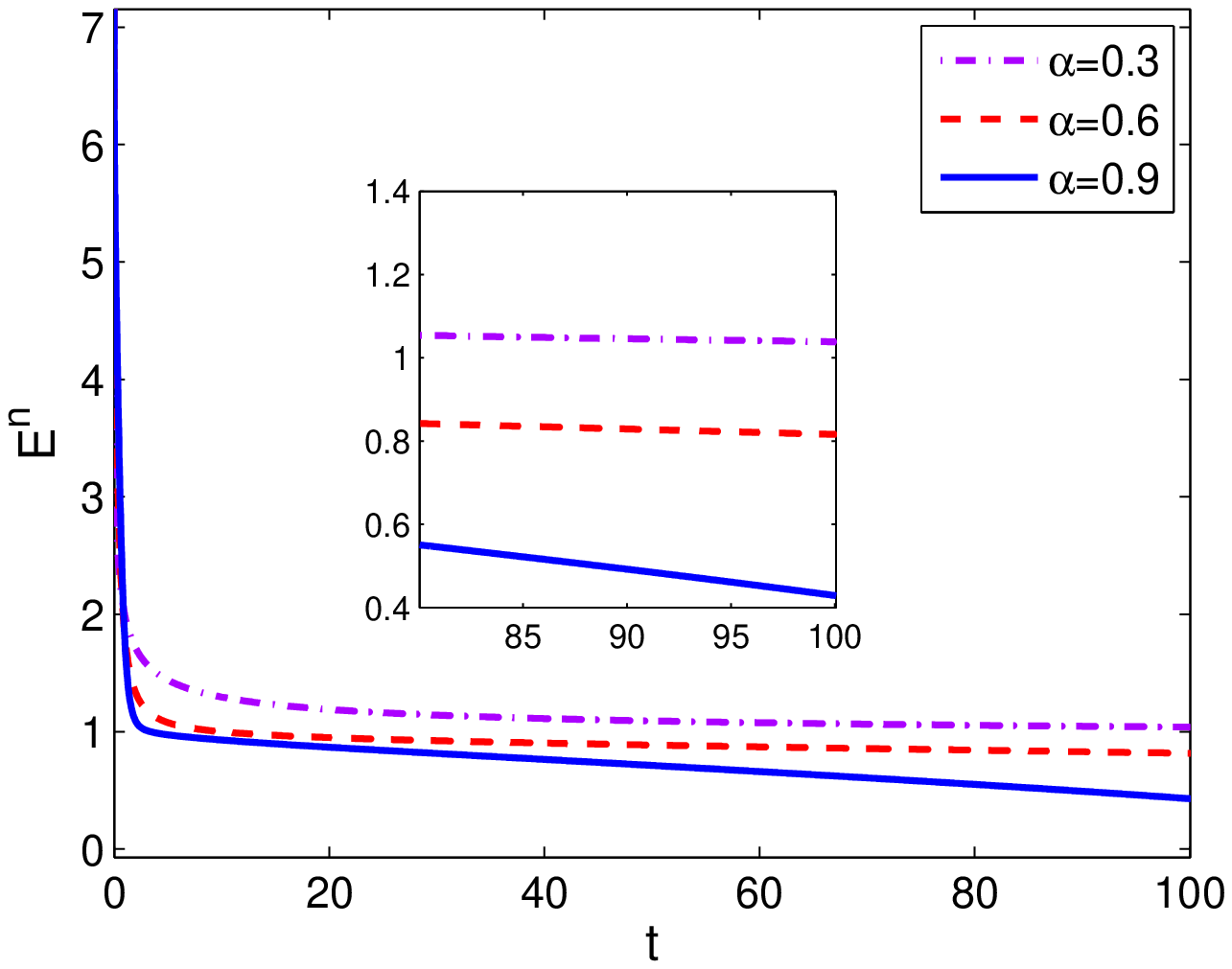}
\includegraphics[width=2.0in]{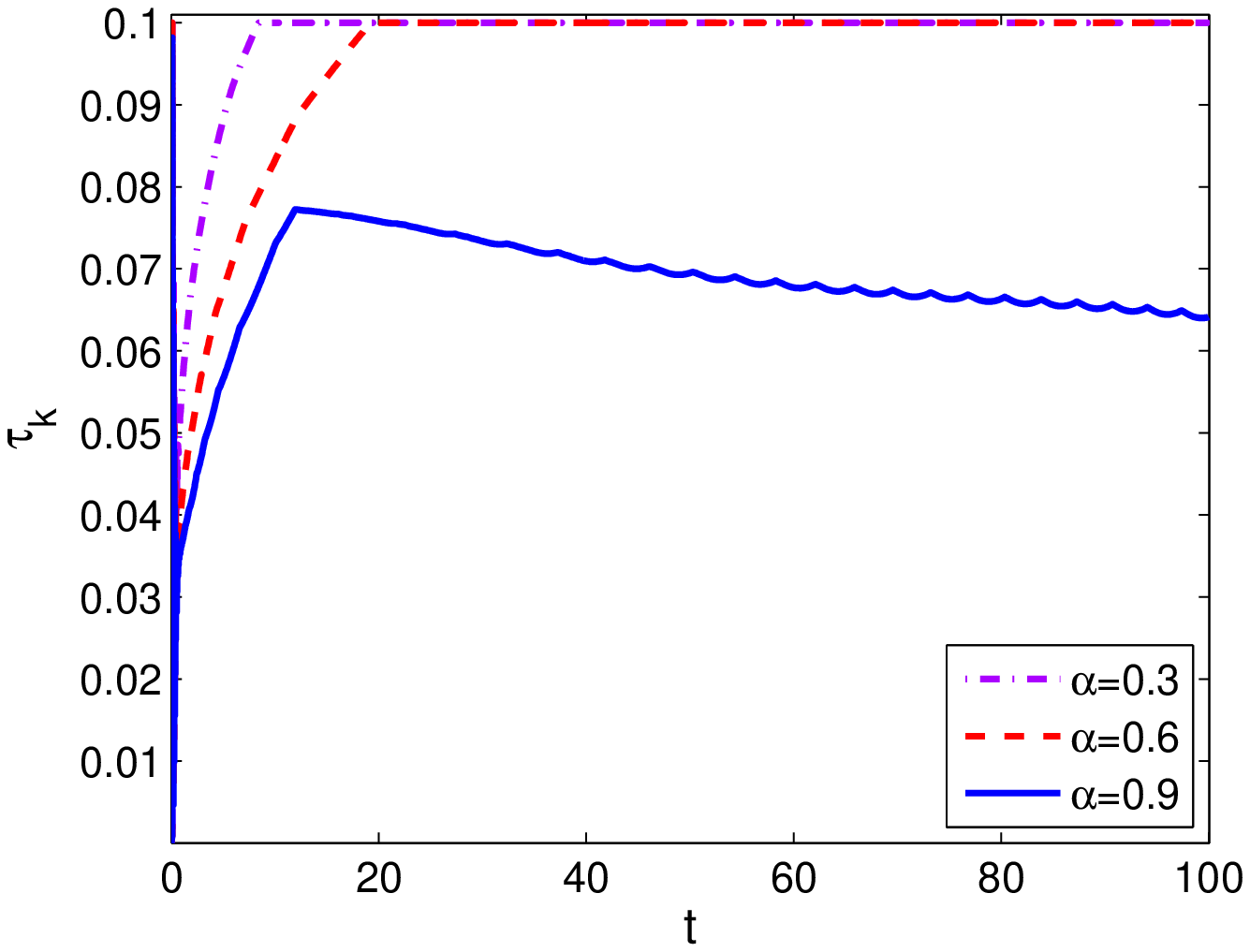}
\caption{The discrete maximum principle (left), energy dissipation (middle) and adaptive time-steps (right) of backward Euler scheme \eqref{L1-1}-\eqref{L1-2}.}
\label{L1MaximumEnergy}
\end{figure}

\begin{figure}[htb!]
\centering
\includegraphics[width=2.0in]{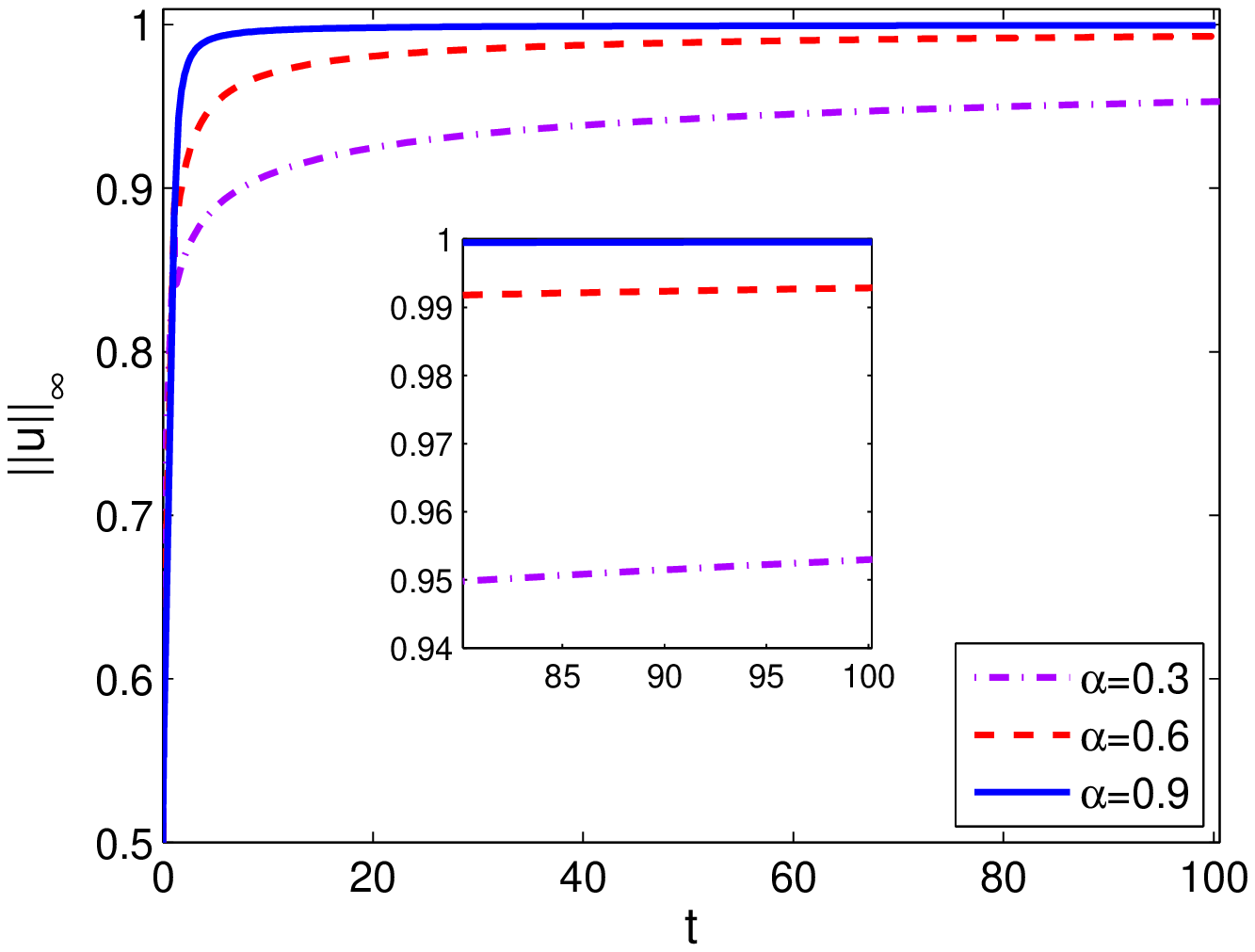}
\includegraphics[width=2.0in]{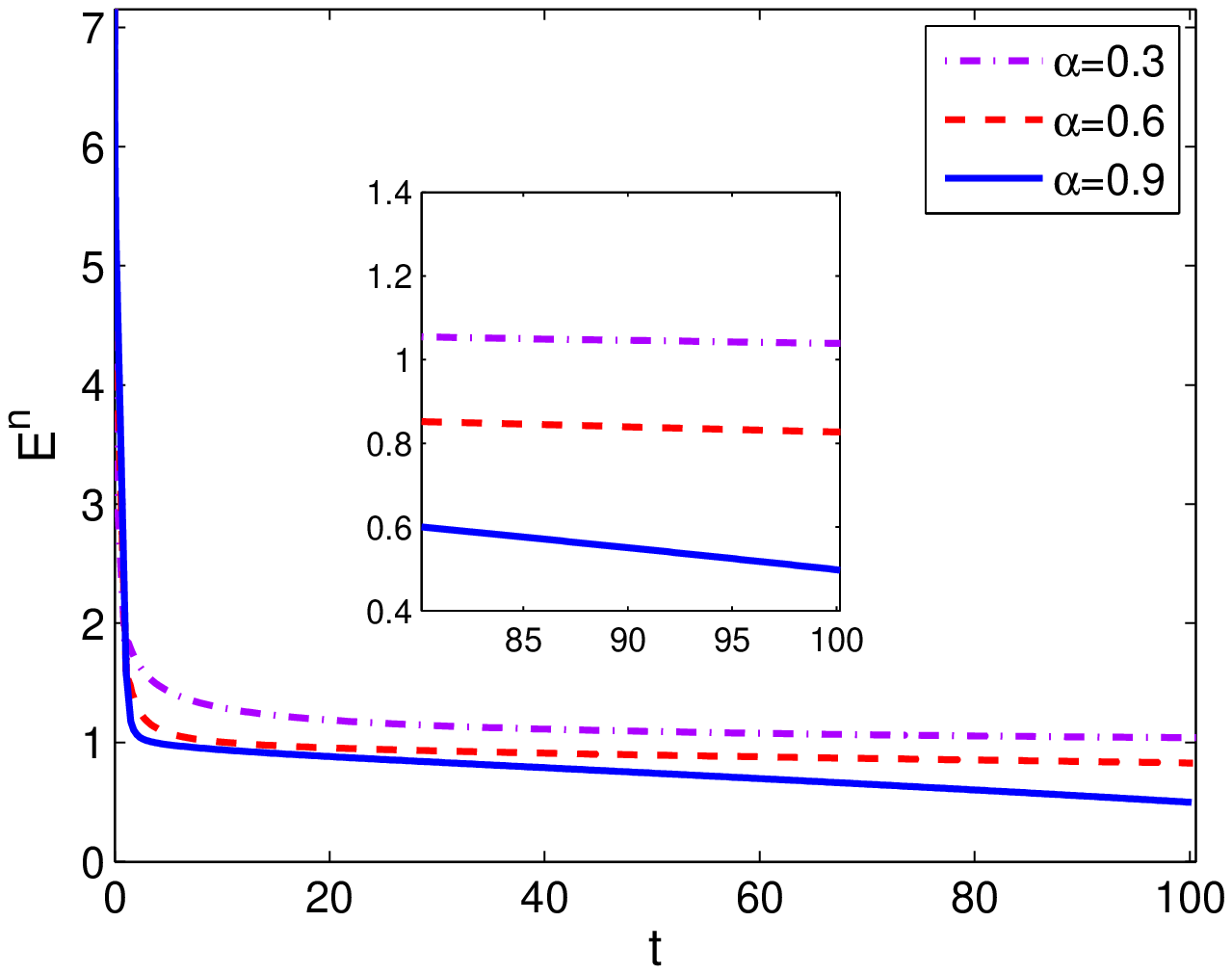}
\includegraphics[width=2.0in]{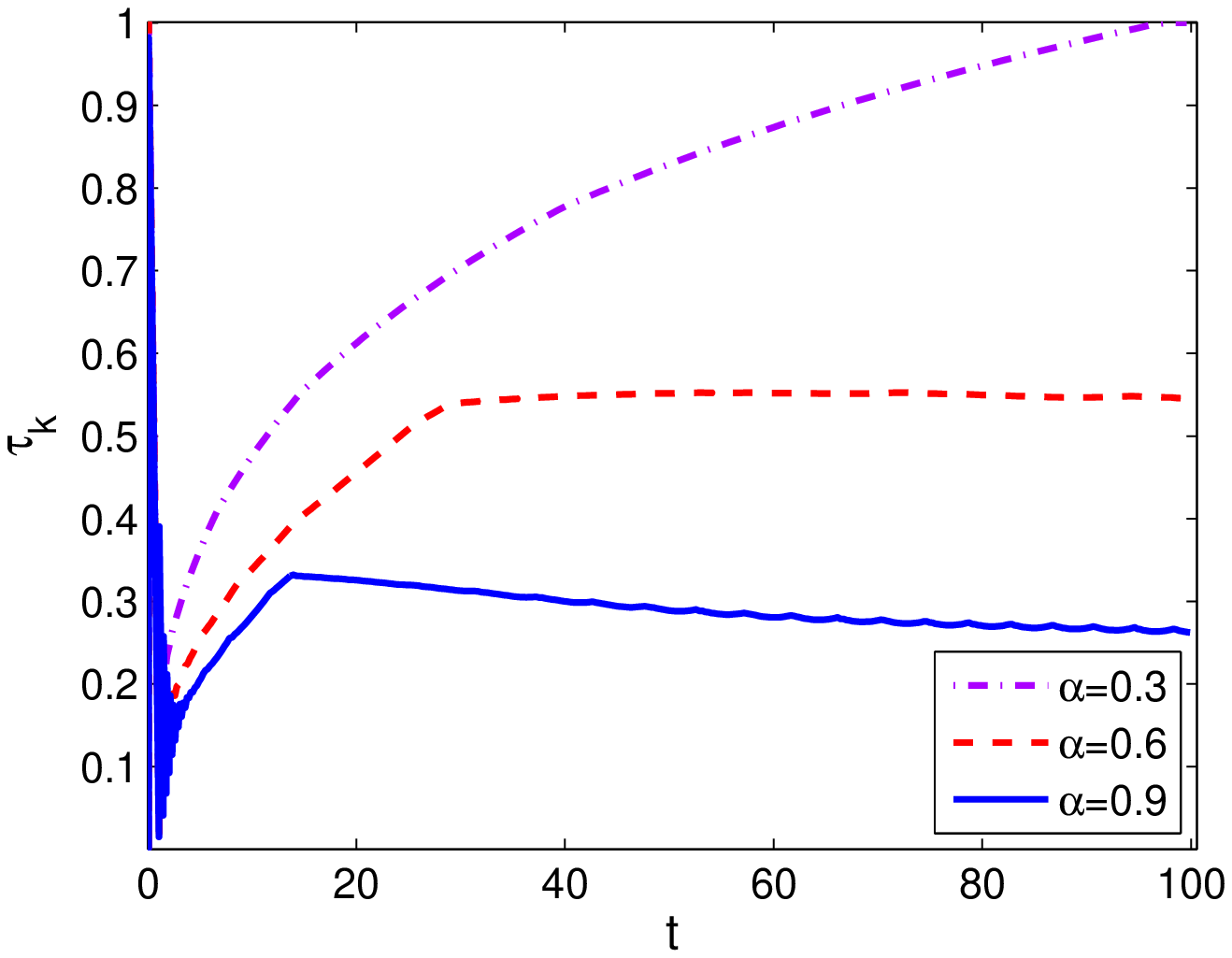}
\caption{The discrete maximum principle (left), energy dissipation (middle) and adaptive time-steps (right)
of stabilized semi-implicit scheme \eqref{StabilizedL1-1}-\eqref{StabilizedL1-2}.}
\label{StabilizedL1MaximumEnergy}
\end{figure}

We find that the solution profiles, generated by the backward Euler scheme \eqref{L1-1}-\eqref{L1-2} with $\tau_{\min}=\tau_{N_{0}}=0.001,\tau_{\max}=0.1,tol=0.15,\beta=200$
and stabilized semi-implicit scheme \eqref{StabilizedL1-1}-\eqref{StabilizedL1-2} with $S=0.1, \tau_{\min}=\tau_{N_{0}}=0.001,\tau_{\max}=1,tol=1.5,\beta=200$ in the remainder interval, are quite identical. Fig. \ref{SnapshotsNumAPPL} gathers some
snapshots at four different times. It is seen that the two bubbles coalesce into a single bubble
as the time escapes, while the rate of coalescence is deeply affected by the fractional order $\alpha$, see \cite{Liu2018Time,Zhao2019On}.
The larger the fractional order $\alpha$, the faster the bubbles coalesce.

Fig.\ref{L1MaximumEnergy} depicts the solution in the maximum norm and the discrete energy ($E^{n}$ is the discrete counterpart of the energy functional defined in the model in spite of no theoretical proof is available in current work)
of the backward Euler scheme \eqref{L1-1}-\eqref{L1-2}.
It is obvious that the solutions are uniformly bounded by the value 1 for different fractional orders $\alpha$,
as predicted by Theorem \ref{DisMaxPrinciple}.
Moreover, the larger the fractional order $\alpha$, the faster it approaches the maximum value.
The middle of Fig. \ref{L1MaximumEnergy} says that the discrete energy
is also decreasing as the time escapes, although we can not verify it theoretically.
The right side of Fig. \ref{L1MaximumEnergy} depicts the adopted time-steps,
 and we observe that the time-steps are always small at the early stage, implying the fast evolution dynamics near the initial time.
Fig.\ref{StabilizedL1MaximumEnergy} shows analogous plots for the stabilized scheme \eqref{StabilizedL1-1}-\eqref{StabilizedL1-2},
where we see the similar behaviors on the maximum norm value, the discrete energy and the adaptive time-steps. Note that the maximum time step $\tau_{\max}=0.1$ of the backward Euler scheme is to ensure the convergence of iterative method, thus we can expect the stabilized scheme to be more efficient than the nonlinear one.

\section{Concluding remarks}
In simulating the time-fractional phase field equations including the Allen-Cahn equation considered in this paper,
the initial singularity should be treated properly because it always destroys the time accuracy of numerical algorithms especially
near the initial time. We consider two fast L1 time-stepping methods on a general class of nonuniform time meshes
such that they will be suitable for both the refined mesh near $t=0$ and certain adaptive time-stepping strategy
to resolve the multiple time scales away from $t=0$.

We show that the nonuniform fast L1 formula can be employed to construct some time-stepping methods preserving the discrete maximum principle
by virtue of the uniform monotonicity of discrete kernels. By using the discrete fractional Gr\"onwall inequality
and global consistency analysis, we established obtain sharp maximum norm error estimates of proposed schemes
and validated them numerically.

It seems challenging to build time-stepping approaches maintaining the discrete energy dissipation law,
especially on general nonuniform time meshes.
Nonetheless, the energy stable schemes permitting adaptive time-stepping strategies are very attractive because
they would be applicable for other time-fractional phase-field models
and for long-time simulations approaching the steady state. These issues will be addressed in the forthcoming reports.

\section*{Acknowledgements}
The authors would like to thank Prof. Jia Zhao and Prof. Yuezheng Gong for their valuable discussions and fruitful suggestions.

\appendix
\section{Discrete fractional Gr\"{o}nwall lemma}\label{sec:gronwall}

The recently developed discrete fractional Gr\"{o}nwall inequality
in \cite{Liao2018discrete} is applicable for any nonuniform time meshes
and suitable for a variety of discrete fractional derivatives.
The following lemma,
involving the Mittag--Leffler function $E_\alpha(z):=\sum_{k=0}^\infty\frac{z^k}{\Gamma(1+k\alpha)}$,
gathers three previous (slightly simplified) results from
\cite[Lemma~2.2, Theorems~3.1 and 3.2]{Liao2018discrete}.
\begin{lemma}\label{FractGronwall}
For $n=1,2,\cdots,N$, assume that the discrete convolution kernels $\{A_{n-k}^{(n)}\}_{k=1}^n$ satisfy the following two assumptions:\\
\textbf{Ass1}. There is a constant $\pi_{a}>0$ such that $A_{n-k}^{(n)}\geq \frac{1}{\pi_{a}}\int_{t_{k-1}}^{t_{k}}\frac{\omega_{1-\alpha}(t_{n}-s)}{\tau_{k}}\zd s$ for $1\leq k\leq n$.\\
\textbf{Ass2}. The discrete kernels are monotone, i.e. $A_{n-k-1}^{(n)}-A_{n-k}^{(n)}\geq 0$ for $1\leq k\leq n-1$.\\
Define also a sequence of discrete complementary convolution kernels~$\{p_{n-j}^{(n)}\}_{j=1}^n$ by
\begin{align}\label{eq: discreteConvolutionKernel-RL}
p_{0}^{(n)}:=\frac{1}{A_{0}^{(n)}},\quad
p_{n-j}^{(n)}:=
\frac{1}{p_0^{(j)}}
\sum_{k=j+1}^{n}\brab{A_{k-j-1}^{(k)}-A_{k-j}^{(k)}}p_{n-k}^{(n)},
    \quad 1\leq j\leq n-1.
\end{align}
Then the discrete  complementary kernels $p^{(n)}_{n-j}\ge0$ are well-defined and fulfill
\begin{align}
&\sum_{j=k}^np^{(n)}_{n-j}A_{j-k}^{(j)}= 1,
\quad\text{for $1\le k\le n\le N$.}\label{eq: P A}\\
&\sum_{j=1}^np^{(n)}_{n-j}\omega_{1+m\alpha-\alpha}(t_j)\leq \pi_{a}\omega_{1+m\alpha}(t_n),
\quad\text{for $m=0,1$ and $1\le n\le N$.}\label{eq: P bound}
\end{align}
Suppose that the offset parameter $0 \leq \nu <1$, $\lambda$ is a non-negative constant independent of the time-steps
and the maximum step size $\tau\le1/\sqrt[\alpha]{2\Gamma(2-\alpha)\lambda\pi_{a}}.$
If the non-negative sequences $(v^k)_{k=0}^N$, $(\xi^{k})_{k=1}^{N}$ and $(\eta^{k})_{k=1}^{N}$ satisfy
\begin{equation}\label{eq: first Gronwall}
\sum_{k=1}^nA_{n-k}^{(n)}\triangledown_{\tau} v^k\le
	\lambda v^{n-\nu}+\xi^n+\eta^n\quad\text{for\ $1\le n\le N$,}
\end{equation}
or
\begin{equation}\label{eq: second Gronwall}
\sum_{k=1}^na^{(n)}_{n-k}\triangledown_{\tau} (v^k)^{2}\le
	\lambda (v^{n-\nu})^{2}+v^{n-\nu}(\xi^n+\eta^n)\quad\text{for\ $1\le n\le N$,}
\end{equation}
then it holds that, for $1\le n\le N$,
\begin{align*}
v^n&\le2E_\alpha\big(2\max\{1,\rho\}\lambda\pi_{a} t_{n}^{\alpha}\big)
	\Big(v^0+\max_{1\le k\le n}\sum_{j=1}^kp^{(k)}_{k-j}(\xi^j+\eta^{j})\Big)\\
&\le2E_\alpha\big(2\max\{1,\rho\}\lambda\pi_{a} t_{n}^{\alpha}\big)
	\Big(v^0+\Gamma(1-\alpha)\pi_{a}\max_{1\le k\le n}\{t_k^{\alpha}\xi^{k}\}+\pi_{a}\omega_{1+\alpha}(t_n)\max_{1\le k\le n}\eta^{k}\Big).
\end{align*}
\end{lemma}


\end{document}